\documentclass[10pt]{cmip-l}
\usepackage[all]{xy}
\usepackage{amsmath}
\usepackage{amsfonts}
\usepackage{amssymb}
\usepackage{amscd}
\usepackage{amsthm}
\usepackage{latexsym}
\usepackage{amsbsy}
\newtheorem{lem}{Lemma}[section]
\newtheorem{prop}{Proposition}[section]
\newtheorem{theorem}{Theorem}[section]

\newtheorem{definition}{Definition}[section]

\newcommand{\A}{\mathcal{A}}

\newcommand{\ot}{\otimes}
\begin{document}
\title{A Short Survey of Cyclic Cohomology}
\author {Masoud Khalkhali}
\address{Mathematics Department,
The University of Western Ontario, 
London, Ontario, N6A 5B7, Canada}

\email{masoud@uwo.ca}

\copyrightinfo{~2010}{Masoud Khalkhali}
\dedicatory{Dedicated with admiration and affection to Alain Connes}

\subjclass[2010]{Primary 58B34; Secondary 19D55, 16T05, 18G30}								

\maketitle

\begin{abstract} This is a short   survey of some aspects of Alain Connes' contributions to  cyclic cohomology theory in the course
of  his work on noncommutative geometry over the  past 30 years. 
\end{abstract}

\tableofcontents

 \section{Introduction}

{\it Cyclic cohomology} was discovered by Alain Connes no later than  1981 and in fact it was
 announced in that year  in a conference in Oberwolfach
\cite{ac81}. I have reproduced the text of his  abstract below. As
it appears in his report, one of Connes'  main motivations to introduce  cyclic cohomology theory  came
from index theory on foliated spaces. Let $(V, \mathcal{F})$ be a
compact foliated manifold and let  $ V/\mathcal{F}$ denote the
space of leaves of $(V, \mathcal{F})$. This space, with its natural quotient topology, is, in general,
 a highly singular
space and in noncommutative geometry one usually replaces the
quotient space $ V/\mathcal{F}$ with a noncommutative algebra $A=
C^* (V, \mathcal{F})$ called the  foliation algebra of $(V,
\mathcal{F})$. It  is the convolution algebra of the holonomy
groupoid of the foliation and is a $C^*$-algebra. It has a dense
subalgebra $\mathcal{A}= C^{\infty} (V, \mathcal{F})$ which plays
the role of the algebra of smooth functions on $ V/\mathcal{F}$.
Let $D$ be a transversally elliptic operator on $(V,
\mathcal{F})$. The analytic index of $D$,  $\text{index} (D) \in
K_0 (A)$,  is an element of the $K$-theory of $A$. This should be compared with the family index theorem \cite{atising} where
the analytic index of a family of fiberwise elliptic operators  is an element of the $K$-theory of the base.
Connes realized
that to identify this class by a cohomological expression it would be necessary  to have a
noncommutative analogue of the Chern character, i.e.,  a map from
$K_0 (\mathcal{A})$ to a,  then  unknown, cohomology theory for
the noncommutative algebra $\mathcal{A}$. This theory, now known as cyclic cohomology,   would then
play
the role of the noncommutative analogue of de Rham homology of  currents for  smooth
manifolds. Its dual version, {\it cyclic homology},  corresponds,
in the commutative case, to de Rham cohomology.

Connes arrived at his definition of cyclic cohomology
 by a  careful analysis of  the algebraic structures
deeply hidden in the (super)traces of products of commutators of operators.
These expressions are directly defined in terms of an elliptic
operator and its parametrix and give  the index of the operator when paired with a
$K$-theory class. In his own words \cite{ac81}:

{\it ``The transverse elliptic theory for foliations requires as a
preliminary step a purely algebraic work, of computing for a
noncommutative algebra $\mathcal{A}$ the cohomology of the following
complex: $n$-cochains are multilinear functions \\$\varphi (f^0, \dots,
f^n)$ of $f^0, \dots ,f^n \in  \mathcal{A}$ where
$$\varphi (f^1, \dots,
f^0)=(-1)^n \varphi (f^0, \dots,
f^n)$$
and the boundary is
\begin{eqnarray*}
b\varphi (f^0, \dots,f^{n+1}) &=&\varphi (f^0f^1, \dots,f^{n+1})-\varphi (f^0, f^1f^2,\dots,
f^{n+1})+\cdots \\& &+(-1)^{n+1} \varphi (f^{n+1}f^0, \dots,
f^n).
\end{eqnarray*}
The basic class associated to a transversally elliptic operator,
for $\mathcal{A}=$ the algebra of the foliation,  is given by:
$$ \varphi (f^0, \dots,f^n)=Trace \, (\varepsilon F[F, f^0][F, f^1]\cdots[F, f^n]), \quad
f^i\in \mathcal{A}$$
where
$$F=\left(
\begin{matrix}
0& Q \\
P & 0
\end{matrix}
\right), \quad  \varepsilon= \left(
\begin{matrix}
1& 0 \\
0& -1
\end{matrix}
\right),
$$
and $Q$ is a parametrix of $P$. An operation
$$S: H^n(\mathcal{A})\to
H^{n+2}(\mathcal{A})$$
is constructed as well as a pairing
$$K(\mathcal{A})
\times H(\mathcal{A})\to \mathbb{C}$$ where $K(\mathcal{A})$ is
the algebraic $K$-theory of $A$. It gives the index of the
operator from its associated class $\varphi$. Moreover $\langle e, \,
\varphi \rangle = \langle e, \, S \varphi \rangle$, so that the important group to determine
is the inductive limit} $H_p=\underset{\rightarrow}{\text{Lim}}\,\it H^n(\mathcal{A})$
{\it for the map $S$. Using the tools of homological algebra the groups
$H^n(\mathcal{A}, \mathcal{A}^*)$ of Hochschild cohomology with
coefficients in the bimodule $\mathcal{A}^*$ are easier to
determine and the solution of the problem is obtained
in two steps:\\
1) the construction of a map
$$B: H^n(\mathcal{A},
\mathcal{A}^*) \to H^{n-1}(\mathcal{A})$$
and the proof of a long exact
sequence
$$
\cdots\to H^n(\mathcal{A}, \mathcal{A}^*) \overset{B}{\to} H^{n-1}(\mathcal{A})\overset{S}{\to}
H^{n+1}(\mathcal{A})\overset{I}{\to}H^{n+1}(\mathcal{A},
\mathcal{A}^*)\to \cdots
$$
where $I$ is the obvious map from the cohomology of the above
complex to
the Hochschild cohomology;\\
2) the construction of a spectral sequence with $E_2$ term given
by the cohomology of the degree $-1$ differential $I\circ B$ on
the Hochschild groups $H^n(\mathcal{A}, \mathcal{A}^*)$ and which
converges strongly to a graded group associated to the inductive
limit.

This purely algebraic theory is then used. For
$\mathcal{A}=C^{\infty}(V)$ one gets the de Rham homology of
currents, and for the pseudo-torus, i.e. the algebra of the
Kronecker foliation, one finds that the Hochschild cohomology
depends on the Diophantine nature of the rotation number while the
above theory gives $H^0_p$ of dimension $2$ and $H^1_p$ of dimension
$2$,  as expected, but from some
remarkable
cancellations."}\\

A full exposition of these results later appeared in two IHES
preprints \cite{ac83a}, and were eventually published as  \cite{ac85}. With the
appearance of \cite{ac85} one could say that the first stage of the
development of noncommutative geometry and specially cyclic
cohomology  reached a stage
of maturity. In the next few sections I shall try to give a quick and concise survey of some aspects
of cyclic cohomology theory as they were developed in \cite{ac85}. The last two
sections are devoted to developments in the subject  after  \cite{ac85}
arising from the work of Connes.

It is a distinct  honor and a great pleasure   to dedicate this short
survey  of cyclic cohomology theory as a small token of our friendship 
to Alain Connes, the originator of the subject,  on the occasion of his  60th birthday.  It inevitably only covers
part of  what  has been done by Alain in this  very important
corner of noncommutative geometry. It is  impossible  to cover everything, and in particular I have left out
many important topics  developed by him including, among others,  the Godbillon-Vey invariant  and type III factors \cite{ac83c}, the
transverse fundamental class for foliations \cite{ac83c}, the Novikov conjecture for hyperbolic groups \cite{como0},
entire cyclic cohomology \cite{ac88},
and multiplicative  characteristic classes \cite{acmk}.  Finally I would
like to thank Farzad Fathi zadeh for carefully reading the text and for
several  useful comments, and Arthur Greenspoon who kindly  edited the whole text.

\section{Cyclic cohomology}
Cyclic cohomology can be  defined in several ways, each shedding
light on a different   aspect of it. Its original definition
\cite{ac81, ac85} was through  a remarkable subcomplex of the
Hochschild complex that we recall first. By algebra in this paper
we mean an associative algebra over $\mathbb{C}$. For an algebra
$\A$ let
$$C^n (\A)=\text{Hom} (\A^{\otimes (n+1)}, \, \mathbb{C}), \quad n=0, 1, \dots,$$
denote the space of $(n+1)$-linear functionals on $\A$. These are our
{\it $n$-cochains}. The Hochschild differential  $b: C^n (\A) \to C^{n+1}
(\A)$ is defined as
\begin{eqnarray*}
(b \varphi)(a_0, \dots ,a_{n+1})&=&\sum_{i=0}^n(-1)^i \varphi (a_0, \dots,a_ia_{i+1},\dots ,a_{n+1})\\
& &+(-1)^{n+1} \varphi (a_{n+1}a_0, \dots, a_{n}).
\end{eqnarray*}
The cohomology of the complex $(C^{*} (\A), b)$ is the
Hochschild cohomology of $\A$ with coefficients in the bimodule $\A^*$. 

The following definition is fundamental and marks the  departure
from Hochschild cohomology in \cite{ac81, ac85}:
\begin{definition}An $n$-cochain $\varphi  \in C^n(\A)$ is called {\it cyclic} if
$$ \varphi (a_n, a_0, \dots ,a_{n-1})=(-1)^n \varphi (a_0, a_1, \dots, a_n)$$
for all $a_0, \dots, a_n$ in $\A$. The space of cyclic
$n$-cochains  will be denoted by  $C^n_{\lambda}(\A).$
\end{definition}
Just why, of all possible symmetry conditions on cochains,   the cyclic property is a reasonable choice
is  at first glance not  at all clear.
\begin{lem} The space of cyclic cochains is invariant under the action
of $b$, i.e.,  
$b \,C^n_{\lambda} (\A) \subset C^{n+1}_{\lambda} (\A)$  for all $n\geq 0.$ 
\end{lem}
To see this one introduces   the operators $\lambda :C^n(\A) \rightarrow
C^n(\A)$ and $b': C^n(\A) \rightarrow C^{n+1}(\A)$ by
\begin{eqnarray*}
(\lambda \varphi)(a_0, \dots ,a_n)&=&(-1)^n \varphi (a_n, a_0, \dots ,a_{n-1}),\\
 (b' \varphi )(a_0, \dots ,a_{n+1})&=&\sum_{i=0}^n(-1)^i \varphi (a_0,
 \dots,a_ia_{i+1},\dots ,a_{n+1}),
\end{eqnarray*}
and  checks that
$ (1-\lambda)b=b'(1-\lambda).$
Since
$ C^{*}_{\lambda} (\A)= \text{Ker}\, (1-\lambda),$
the lemma is proved.

 We therefore have a subcomplex of the Hochschild complex, called the
{\it cyclic complex} of $\A$:
\begin{equation} \label{cc34} C_{\lambda}^0 (\A)\overset{b}{\longrightarrow} C_{\lambda}^1 (\A)
\overset{b}{\longrightarrow}C^2_{\lambda}(\A) \overset{b}{\longrightarrow}\cdots .
\end{equation}

\begin{definition} The cohomology of the cyclic complex \eqref{cc34} is the {\it cyclic cohomology} of
$\A$ and will be denoted by $HC^n (\A)$, $n=0, 1, 2,\dots.$
\end{definition}

And that is Connes' first definition of cyclic cohomology. A cocycle for the cyclic cohomology group $HC^n (\A)$ is called a {\it
cyclic n-cocycle} on $\A$. It is an $(n+1)$-linear functional $\varphi $ on $\A$
which
satisfies the two conditions:
$$ (1-\lambda )\varphi =0, \quad \text{and} \quad b \varphi=0.$$

The inclusion of complexes
\begin{equation} \label{cytoho}(C^{*}_{\lambda}(\A), \, b) \hookrightarrow (C^{*}(\A), \, b)
\end{equation}
induces a map $I$ from   cyclic cohomology  to
Hochschild cohomology:
$$ I: HC^n(\A) \longrightarrow HH^n(\A), \quad n=0, 1, 2, \dots.$$
A  closer inspection of the long exact sequence associated
to \eqref{cytoho},  yields {\it Connes'  long exact sequence}
relating Hochschild cohomology to cyclic cohomology. This is
however easier said than done. The reason is that to identify the
cohomology of the quotient one must use another long exact
sequence, and combine the two long exact sequences to obtain the
result. To simplify the notation, let us denote the Hochschild and
cyclic complexes by $C$ and  $C_{\lambda}$, respectively. Then
\eqref{cytoho} gives us an  exact
 sequence of complexes
\begin{equation}\label{ses0} 0 \to C_{\lambda}\to C \overset{\pi}{\to} C
/C_{\lambda} \to 0.
\end{equation}
Its associated long exact sequence is
\begin{equation} \label{ses1} \cdots \longrightarrow HC^n (\A) \longrightarrow HH^n (\A)
\longrightarrow H^n (C/C_{\lambda})
\longrightarrow  HC^{n+1}(\A) \longrightarrow \cdots
\end{equation}
We need to identify the cohomology groups $H^n(C/C_{\lambda})$. To this
end, consider  the short exact sequence of complexes
\begin{equation}\label{ses2} 0 \longrightarrow  C/C_{\lambda} \overset{1-\lambda}{\longrightarrow} (C, b')
\overset{N}{\longrightarrow}
C_{\lambda} \longrightarrow
 0 ,
 \end{equation}
 where the  operator $N$ is defined by
 $$ N=1+\lambda +\lambda^2 +\cdots +\lambda^n : C^n \longrightarrow
 C^n.$$  The relations $ (1-\lambda)b=b'(1-\lambda),$
 $ N (1-\lambda)=(1-\lambda) N =0,$  and   $bN=Nb'$
  show that $1-\lambda$ and $N$ are morphisms of
 complexes in \eqref{ses2}. As for the exactness of \eqref{ses2}, the
 only nontrivial part is to show that
 $\text{ker}\, (N) \subset \text{im}\, (1-\lambda)$, which can be verified.
  Assuming $\A$ is unital,  the middle complex $(C, \, b')$  in  \eqref{ses2}
   can be shown to be  exact with
   a contracting
 homotopy
 $ s: C^n \to C^{n-1}$
 defined  by
 $ (s \varphi) (a_0, \dots, a_{n-1})=(-1)^{n-1} \varphi (a_0, \dots,
 a_{n-1}, 1),
 $
which satisfies
$b's +sb' =\text{id}.$
The long exact sequence associated to \eqref{ses2}  looks like
\begin{equation} \label{ses3} \cdots \longrightarrow H^n (C/C_{\lambda}) \longrightarrow H^n_{b'} (C)
\longrightarrow HC^n (\A)
\longrightarrow  H^{n+1}(C/C_{\lambda}) \longrightarrow H^{n+1}_{b'} (C)\longrightarrow \cdots
\end{equation}
Since $H^n_{b'} (C)=0$ for all $n$, it follows that the connecting homomorphism
 \begin{equation} \label{iso21} \delta: HC^{n-1} (\A) \to  H^n(C/C_{\lambda})
 \end{equation}
  is an {\it isomorphism} for all $n\geq 0.$  Using this in \eqref{ses1},  we obtain {\it Connes'
   long exact sequence} relating Hochschild  and cyclic cohomology:
  \begin{equation}\label{ibs}
   \cdots \longrightarrow  HC^n(\A) \overset{I}{\longrightarrow}HH^n(\A) \overset{B}{\longrightarrow}
  HC^{n-1}(\A) \overset{S}{\longrightarrow} HC^{n+1}(\A) \longrightarrow
  \cdots .
  \end{equation}

The operators $B$ and $S$ play a prominent role in
 noncommutative geometry.
As we shall see, the operator $B$ is the
analogue of de Rham's differential in the noncommutative world,
while the  {\it periodicity operator}  $S$  is closely related to Bott
periodicity in topological $K$-theory.  Remarkably, there is a formula for $B$ on the
level of cochains given by
$ B = NB_0,$
where $B_0: C^n \to C^{n-1}$ is defined by
$$ B_0 \varphi (a_0, \dots, a_{n-1})=\varphi (1, a_0, \dots, a_{n-1})
-(-1)^n \varphi (a_0, \dots, a_{n-1}, 1).$$
Using the relations $ (1-\lambda) b=b' (1-\lambda), \, (1-\lambda)N= N(1-\lambda)=0,\,  b N= Nb',$
and $sb'+b's=1$, one shows  that
\begin{equation}\label{bB}bB+ Bb=0, \quad \text{and} \quad B^2=0.
\end{equation}

  Using the periodicity operator $S$, the {\it periodic
cyclic cohomology} of $\A$ is then  defined as the  direct limit of cyclic cohomology groups under the
operator $S$:
$$ HP^i (\A): = \underset{\longrightarrow}{\text{Lim}} \, HC^{2n+i} (\A), \quad \quad i=0, \, 1.$$
Notice that since $S$ has degree 2, there are only two   periodic
groups. These periodic groups have better stability properties
compared to cyclic cohomology groups. For example, they are homotopy
invariant, and they pair with K-theory.

 A much deeper relationship  between Hochschild and cyclic
cohomology groups is encoded in Connes' {\it $(b, B)$-bicomplex} and the associated
Connes spectral sequence that we shall briefly
recall now. Consider the relations \eqref{bB}.
The \index{$(b, B)$-bicomplex} $(b, B)$-bicomplex of a unital algebra
 $\A$, denoted by $\mathcal{B} \, (\A)$,   is the bicomplex
$$\begin{CD}
\vdots @.\vdots @.\vdots \\
C^2(\A)@>B>> C^1(\A)@>B>> C^0(\A)\\
@AbAA @AbAA\\
C^1(\A)@>B>>C^0(\A)\\
@AbAA\\
C^0(\A)
\end{CD}
$$

As usual, there are two spectral sequences attached to this bicomplex.
The following fundamental result of Connes \cite{ac85} shows that  the spectral sequence
obtained from  filtration by rows converges to cyclic cohomology. Notice that the
$E^1$ term of this
spectral sequence is the
 Hochschild cohomology of $\A$.
 \begin{theorem}  \label{bBbicom}The map
$ \varphi \mapsto  (0, \dots, 0,\varphi)$ is a quasi-isomorphism
of complexes
\rm{$$(C_{\lambda}^{*} (\A),  \, b) \to (\text{Tot}\,  \mathcal{B} \, (\A),
\, b+B).$$}
\end{theorem}
This is a consequence of the  vanishing  of the $E^2$ term  of the
second spectral sequence (filtration by columns) of $\mathcal{B}
(\A)$. To prove this, Connes  considers the short exact sequence of
$b$-complexes
$$ 0 \longrightarrow \text{Im}\, B \longrightarrow \text{Ker}\,  B
\longrightarrow  \text{Ker}\,
B/ \text{Im}\, B \longrightarrow 0, $$
and proves that  (\cite{ac85},  Lemma 41), the induced map
$$ H_b (\text{Im}\, B) \longrightarrow  H_b (\text{Ker}\, B)$$
is an isomorphism. This is a very technical result. It follows that $H_b (\text{Ker}\, B / \text{Im}\, B)$ vanishes. To
take care of the first column one appeals to the fact that
$ \text{Im}\, B \simeq \text{Ker} \, (1-\lambda)$
is the space of cyclic cochains.

We give an alternative proof of Theorem \eqref{bBbicom} above. To this end,
 consider the  {\it cyclic bicomplex} $\mathcal{C} (\A)$  defined by
 $$\begin{CD}
\vdots @.\vdots @.\vdots @.\\
C^2 (\A)@ >1-\lambda >> C^2 (\A) @ > N >> C^2 (\A)@ >1- \lambda >> \cdots  \\
@AAbA  @ AA- b'A @AA  b A  \\
C^1 (\A)@ >1-\lambda >> C^1 (\A) @ > N >> C^1 (\A)@ >1-\lambda >> \cdots \\
@AA  b A @AA- b'A @AA  bA  \\
C^0 (\A) @>1-\lambda >> C^0 (\A) @ >  N >> C^0 (\A)@ >1-\lambda >> \cdots
\end{CD}
$$
The  total cohomology of $\mathcal{C} (\A)$ is
isomorphic to cyclic cohomology:
$$ H^n (\text{Tot} \,\mathcal{C}\, (\A))\simeq HC^n (\A), \quad n\geq 0.$$
This is a consequence  of the simple fact that the rows of
$\mathcal{C} (\A)$ are exact except in
degree zero, where  their cohomology  coincides with the
 cyclic complex $(C^{*}_{\lambda}(\A), b)$.
 So it suffices to show that $\text{Tot}\, \mathcal{B}(\A)$ and $ \text{Tot}\,\mathcal{C}(\A)$
are quasi-isomorphic. This can be done  by  explicit
formulas. Consider the chain maps 
\begin{eqnarray*}
I: \text{Tot}\, \mathcal{B}(\A) &\rightarrow &  \text{Tot}\,\mathcal{C}(\A), \quad I=\text{id} +Ns,\\
J: \text{Tot}\, \mathcal{C}(\A) &\rightarrow & \text{Tot}\,
\mathcal{B}(\A), \quad  J=\text{id}
+sN.
\end{eqnarray*}
It can be directly  verified that the following operators  define chain
homotopy equivalences:
\begin{eqnarray*}
g: \text{Tot}\, \mathcal{B}(A) &\rightarrow & \text{Tot}\, \mathcal{B}(A),  \quad g=Ns^2B_0,\\
h: \text{Tot}\, \mathcal{C}(A) & \rightarrow & \text{Tot}\,
\mathcal{C}(\A), \quad h=s.
\end{eqnarray*}

To give an example of an application of the spectral sequence of Theorem \eqref{bBbicom}, let me
 recall Connes' computation of the {\it continuous cyclic cohomology}  of the topological algebra
   $\A=C^{\infty}(M)$, i.e.,    the algebra of smooth
complex valued
functions on a closed smooth  $n$-dimensional  manifold $M$. This
example is important since, apart from its applications,  it clearly
shows that cyclic cohomology is a noncommutative analogue of de Rham
homology.

The continuous analogues of Hochschild and cyclic cohomology for
topological algebras are defined as follows \cite{ac85}. Let $\A$ be a
topological algebra. A {\it continuous cochain} on $\A$ is a jointly continuous multilinear map
$\varphi: \A \times \A \times  \cdots \times \A \to \mathbb{C}$. By working with
just continuous cochains,  as opposed to all cochains, one obtains the
continuous analogues of Hochschild and cyclic cohomology groups. In
working with algebras of smooth functions (both in the commutative and
noncommutative case), it is essential to use this continuous analogue.

The topology of
$C^{\infty}(M)$ is defined by the sequence of seminorms
$$\|f\|_n = \text{sup} \, | \partial^{\alpha} \, f|,  \quad |\alpha | \leq n,$$
where the supremum is over a fixed, finite,  coordinate cover for
$M$. Under this topology, $C^{\infty}(M)$  is a locally convex, in fact nuclear,  topological algebra. Similarly one
topologizes the space of $p$-forms on $M$ for all $p\geq 0$.
Let $$\Omega_p M: =\text{Hom}_{\text{cont}} (\Omega^p M, \mathbb{C})$$
 denote the continuous dual of the space of $p-$forms on $M$. Elements
of   $\Omega_p M$ are called {\it de Rham $p$-currents}.
By dualizing the de Rham differential $d$,
we obtain a differential $ d^* : \Omega_{*} M \to \Omega_{*
-1} M$,  and a complex, called the
{\it de Rham complex of currents} on $M$:
$$ \Omega_0 M \overset{d^*}{\longleftarrow}\Omega_1 M
\overset{d^*}{\longleftarrow}\Omega_2 M \overset{d^*}{\longleftarrow} \cdots . $$
The homology of this complex is  the {\it de Rham homology}
of $M$ and we
 denote it by $H_{*}^{dR} (M)$.

It is easy to check that for any $m$-current $C$, closed or not, the cochain 
\begin{equation} \label{drthoch}
\varphi_C (f_0, f_1, \dots , f_m): =\langle C, \,  f_0df_1 \cdots df_m \rangle,
\end{equation}
is a continuous Hochschild cocycle on $ C^{\infty} (M)$. Now  if $C$ is
closed,   then one  checks that $\varphi_C$ is a cyclic
$m$-cocycle on $C^{\infty} (M)$. Thus we obtain   natural maps
\begin{equation}\label{deHoch} \Omega_m M \to HH^m_{cont} (C^{\infty} (M))
\end{equation}
and
\begin{equation}\label{cycHoch}
  Z_m M \to HC^m_{cont}(C^{\infty} (M)),
\end{equation}
where $Z_m (M)\subset \Omega_m M $ is the space of closed $m$-currents
on $M$.
For example,  if $M$ is oriented and $C$ represents its orientation class, then
\begin{eqnarray} \label{cycle1}\varphi_C (f_0, f_1, \dots , f_n)
=\int_M f_0df_1 \cdots df_n,
\end{eqnarray}
which is easily checked to be  a cyclic $n$-cocycle on $\A$.

 In \cite{ac85}, using an
 explicit resolution, Connes shows that \eqref{deHoch} is a quasi-isomorph-\\ism.  Thus  we have a natural  isomorphism between  space of de
 Rham currents on $M$ and the continuous Hochschild cohomology of $C^{\infty}(M):$
 \begin{equation} \label{ccdrcur1} HH^i_{\text{cont}} (C^{\infty}(M)) \simeq \Omega_i M   \quad \quad i=0, 1,
 \dots
 \end{equation}
To compute the continuous  cyclic homology of $\A$, one first observes that
 under the isomorphism \eqref{ccdrcur1} the
operator $B$ corresponds to the de Rham differential $d^*$. More
precisely, for each integer $n\geq 0$ there is a commutative diagram:
$$\begin{CD}
\Omega_{n+1} M@> \mu >> C^{n+1}(\A) \\
@VVd^*V @VVBV\\
\Omega_{n} M@> \mu >> C^{n}(\A)\\
\end{CD}
$$
where $\mu (C)= \varphi_C$ and $\varphi_C$ is defined by \eqref{drthoch}. Then, using the
spectral sequence of Theorem \eqref{bBbicom} and the isomorphism
\eqref{ccdrcur1}, 
 Connes obtains \cite{ac85}:
\begin{equation}\label{cothm1}
HC^n_{\text{cont}}(C^{\infty} (M))\simeq Z_n (M) \oplus H_{n-2}^{\text{dR}}(M) \oplus \cdots \oplus H_k^{\text{dR}}
 (M),
 \end{equation}
where  $k=0$ if $n$ is even and $k=1$ is $n$ is odd.
For the continuous periodic cyclic cohomology he  obtains
 \begin{equation} \label{cothm2}
 HP^k _{\text{cont}} (C^{\infty} (M))\simeq \bigoplus_i H_{2i+k}^{\text{dR}}(M), \quad k=0, 1.
  \end{equation}

We shall also briefly recall Connes' computation of the Hochschild and cyclic
 cohomology  of smooth noncommutative tori \cite{ac85}. This result
 already appeared in Connes' Oberwolfach report \cite{ac81}.
 When $\theta$ is rational,
 the smooth noncommutative torus  $\mathcal{A}_{\theta}$ can be shown to be  Morita equivalent to
 $C^{\infty} (T^2)$, the algebra of smooth functions on the 2-torus. One can then use 
  Morita invariance of Hochschild and cyclic cohomology to reduce the
  computation of these groups to those for the algebra $C^{\infty}
  (T^2)$. This takes care of  rational $\theta$. So we can
  assume $\theta$ is irrational and we denote the  generators of
  $\mathcal{A}_{\theta}$ by $U$ and $V$ with the relation $V U =\lambda UV$, where
  $\lambda = e^{2\pi i \theta}$.

 Recall that an irrational number $\theta$ is  said to
satisfy a Diophantine condition if $ |1-\lambda^n|^{-1}= O (n^k)$
for some positive integer $k$.
\begin{prop} (\cite{ac85}) \label{ncthh2}Let $\theta \notin \mathbb{Q}$. Then\\
a) One has $HH^0 (\mathcal{A}_{\theta}) =\mathbb{C},$\\
b) If $\theta$ satisfies a Diophantine condition then
$HH^i(\mathcal{A}_{\theta})$ is 2-dimensional for i=1  and  is 1-dimensional for $i=2,$\\
c) If $\theta$ does not satisfy a Diophantine condition, then
$HH^i(\mathcal{A}_{\theta})$ are infinite dimensional non-Hausdorff
spaces for $i=1, 2.$
\end{prop}

 Remarkably, for all values of $\theta$, the periodic cyclic cohomology is finite dimensional and is
 given by
$$ HP^0 (\mathcal{A}_{\theta})= \mathbb{C}^2, \quad \quad  HP^1 (\mathcal{A}_{\theta})= \mathbb{C}^2.$$
An explicit  basis for these groups are given by cyclic 1-cocycles
$$ \varphi_1 (a_0, a_1) = \tau (a_0 \delta_1 (a_1)),  \quad \text{and} \quad \varphi_1
(a_0, a_1) = \tau (a_0 \delta_2 (a_1))$$
and cyclic 2-cocycles
$$ \varphi (a_0, a_1, a_2)=\tau (a_0 (\delta_1(a_1)
\delta_2(a_2)-\delta_2(a_1) \delta_1(a_2))), \quad \text{and} \quad S\tau,$$
where $\delta_1, \delta_2 : \mathcal{A}_{\theta} \to
\mathcal{A}_{\theta}$ are the canonical derivations defined by
$$ \delta_1 (\sum a_{mn}U^mV^n)=\sum ma_{mn} U^mV^n,  \quad \quad \delta_2
(U^mV^n)=\sum n a_{mn} U^mV^n,$$
and $\tau : \mathcal{A}_{\theta} \to \mathbb{C}$ is the canonical trace.

A noncommutative generalization of formulas like \eqref{cycle1}   was
introduced in \cite{ac85} and 
played an important role in the development of cyclic cohomology theory in general.
 It gives a  geometric meaning to the notion of a cyclic cocycle over an algebra and goes as follows.
Let $(\Omega,\, d)$ be a differential graded algebra.
 A {\it closed graded
trace} of dimension $n$ on  $(\Omega, \, d)$ is a linear map
$$\int : \Omega^n \longrightarrow  \mathbb{C}$$
such that
$$ \int d\omega =0, \quad \text{and} \int [\omega_1, \,
\omega_2]=0,$$
for all $\omega$ in $\Omega^{n-1}$, $\omega_1$ in $\Omega^i$,  $\omega_2$ in $\Omega^j$ and $i+j=n.$
 An $n$ dimensional cycle over an algebra $\A$ is a triple  $( \Omega, \, \int, \,\rho)$, where $\int$
 is an
$n$-dimensional  closed  graded trace on $(\Omega, \, d)$ and
$\rho: A \to \Omega_0$ is an algebra homomorphism.
Given a cycle   $(\Omega, \, \int, \,  \rho)$ over $\A$,     its {\it
character} is the
cyclic $n$-cocycle on $\A$ defined by
\begin{equation} \label{chacycle}
\varphi (a_0, a_1, \dots ,a_n)=\int \rho(a_0)d \rho(a_1) \cdots d \rho(a_n).
\end{equation}
Conversely one shows that all cyclic cocycles  are obtained in this way.

Once one has the  definition of cyclic cohomology, it is not
difficult to formulate a dual notion of {\it cyclic homology} and
a pairing between the two.  Let $C_n (\A)= \A^{\otimes (n+1)}$.
The analogues of the operators $b, b'$ and $\lambda$ are easily
defined on $C_*(\A)$ and are usually denoted by the same letters, as we do here.
 For example $b: C_n (\A) \to C_{n-1}(\A)$ is defined by
\begin{eqnarray}\label{hhomb1}
b (a_0\otimes  \cdots \otimes a_n)&=&
 \sum_{i=0}^{n-1} (-1)^i(a_0\otimes  \cdots
\otimes a_ia_{i+1} \otimes \cdots
\otimes a_n)\\
& +&(-1)^{n}(a_n a_0\otimes \cdots
\otimes a_{n-1}).
\end{eqnarray}
Let
$$ C_n^{\lambda}(\A):= C_n(\A)/\text{Im} (1-\lambda).$$
The relation
$ (1-\lambda )b'=b(1-\lambda)$ shows that the operator $b$ is
well defined on
$C_{*}^{\lambda}(\A)$. The  complex
$(C_{*}^{\lambda}(\A), \, b)$
is called the
 {\it homological cyclic complex} of $\A$ and its  homology, denoted by $HC_n (\A), n=0, 1, \dots$,
 is
 the  {\it cyclic homology} of $\A$. The evaluation map
 $\langle \varphi, (a_0 \otimes \dots \otimes a_n)\rangle \mapsto \varphi (a_0, \dots, a_n)$
 clearly defines
 a degree zero pairing $HC^{*}(\A)\otimes HC_{*}(\A) \to \mathbb{C}$.
Many results of cyclic cohomology theory,  such as Connes' long exact sequence   and spectral
 sequence, and
  Morita invariance, continue to hold for cyclic homology theory with
  basically the same proofs.

Another important idea of Connes in the 1980's was the introduction of
{\it entire cyclic cohomology} of Banach algebras \cite{ac88}. This allows one to deal with algebras
of functions on infinite dimensional (noncommutative) spaces such as
those appearing in constructive quantum field theory. These algebras typically don't carry
finitely summable   Fredholm modules, but in some cases have so-called
$\theta$-summable  Fredholm modules. In \cite{ac88} Connes extends the definition of Chern character
to such  Fredholm modules  with values in entire cyclic cohomology. 

After the appearance of \cite{ac85},  cyclic (co)homology  theory took on many lives and was further developed
along distinct lines, including a purely algebraic one,  with a big impact on algebraic $K$-theory. The cyclic
cohomology of many algebras was later computed including the
very important case of group algebras \cite{bur} and groupoid algebras. For many of these more
algebraic aspects of the theory we refer to \cite{lod, kh} and references therein.

\section{From $K$-homology to cyclic cohomology}

As I  said in the introduction, Connes' original motivation  for the
development
of  cyclic cohomology was to give a receptacle for  a  noncommutative
Chern character map
on the {\it $K$-homology} of noncommutative algebras. The cycles of
$K$-homology  can be  represented by,  even or odd,  {\it Fredholm
modules}. Here we just focus on the odd case, and we refer to
\cite{ac85, acbook} for the even case. Given a Hilbert space $\mathcal{H}$,
let $\mathcal{L} (\mathcal{H})$ denote the algebra of bounded
linear operators on $\mathcal{H}$, and $\mathcal{K}(\mathcal{H})$ denote
the algebra of compact operators. Also, for $1\leq p<\infty$, let
$\mathcal{L}^p (\mathcal{H})$ denote the Schatten ideal of
$p$-summable operators. By definition,  $T\in \mathcal{L}^p (\mathcal{H})$
if  $|T|^p$ is a trace class operator.
\begin{definition} \label{fm1} An  odd  Fredholm module over a unital  algebra
$\A $ is a pair  $(\mathcal{H}, \, F)$ where \\
1. $\mathcal{H}$ is a Hilbert space endowed with a representation
$$ \pi : \A \longrightarrow \mathcal{L} (\mathcal{H}),$$
2. $F \in \mathcal{L}(\mathcal{H})$ is a bounded selfadjoint  operator with $F^2=I$, \\
3. For   all $a
\in \A$ we have
\begin{equation} \label{psabl}
[F, \,  \pi(a)]=F\pi(a) -\pi(a)F \in \mathcal{K}(\mathcal{H}).
\end{equation}
\end{definition}

A Fredholm module $(\mathcal{H}, \, F)$ is called  {\it $p$-summable}  if,
instead of \eqref{psabl},   we have the stronger condition:
\begin{eqnarray} \label{fs11}[F, \, \, \pi(a)] \in \mathcal{L}^p (\mathcal{H})
\end{eqnarray}
for all $a \in \A$.

To give a simple example, let $A=C (S^1)$ be the algebra of continuous functions on the circle and let $A$   act on
$\mathcal{H}=L^2 (S^1)$ as multiplication operators.
  Let $F(e_n)=e_n$ if  $n\geq 0$ and $F(e_n)=-e_n$ for $n < 0$,  where
   $e_n (x)= e^{2\pi inx}, n\in \mathbb{Z},$ denotes the standard orthonormal basis of $\mathcal{H}$.
   Clearly $F$ is  selfadjoint
   and $F^2=I$. To show that $[F, \, \pi (f)]$ is a compact operator for
  all $f \in C(S^1)$, notice that if $f = \sum_{|n|\leq N} a_n e_n$ is a
  finite trigonometric sum then $[F,\, \pi ( f)]$ is a finite rank operator and
  hence is compact. In general we can uniformly approximate a continuous function by a
  trigonometric sum and show that
  the commutator is compact for any continuous $f$. This shows that $(\mathcal{H},
  \, F)$ is an odd Fredholm module over $C(S^1)$.
This Fredholm module is not $p$-summable for any $1\leq p< \infty$. If we restrict it
to the subalgebra $C^{\infty} (S^1)$ of smooth functions, then it can be checked that
$(\mathcal{H}, \, F)$ is in fact $p$-summable for all $p >1$, but  is not 1-summable even in this case.

Now let me describe Connes' noncommutative Chern character from
$K$-homology to cyclic cohomology. 
Let $(\mathcal{H}, \, F)$ be an odd $p$-summable Fredholm module over an algebra $\A$. For any odd integer $2n-1$ such
that $2n \geq p$,
Connes defines a cyclic $(2n-1)$-cocycle $\varphi_{2n-1}$ on $\A$ by \cite{ac85}
\begin{equation} \label{oddcc}\varphi_{2n-1} (a_0, a_1, \dots, a_{2n-1})= \text{Tr} \,(F [F,
\,a_0][F,\,
a_1]\cdots [F, \,a_{2n-1}]),
\end{equation}
where Tr denotes the operator trace and instead of $\pi (a)$ we
simply write $a$. Notice that by our  $p$-summability assumption,
each commutator is in $\mathcal{L}^p (\mathcal{H})$  and hence, by
H\"{o}lder inequality for Schatten class operators,  their product
is in fact a trace class operator as soon as $2n\geq p$. One checks by a
direct computation that $\varphi_{2n-1}$ is a cyclic cocycle.

The next proposition shows that these cyclic cocycles are related to each other via the periodicity
$S$-operator of cyclic cohomology.  This is probably how Connes came
across the periodicity operator $S$  in the first place.

\begin{prop} For all $n$ with $2n\geq p$ we have
$$ S \varphi_{2n-1}= -(n+\tfrac{1}{2})\, \varphi_{2n+1}.$$
\end{prop}

By rescaling $\varphi_{2n-1}$'s, one obtains a well defined element in the periodic cyclic
cohomology.   The (unstable) odd {\it Connes-Chern character}
$\text{Ch}^{2n-1}= \text{Ch}^{2n-1} (\mathcal{H}, \,F)$ of an odd finitely
summable Fredholm module $(\mathcal{H}, \,F)$ over $\A$ 
is   defined by  rescaling the cocycles $\varphi_{2n-1}$ appropriately.  Let
$$\text{Ch}^{2n-1} \,(a_0, \dots, a_{2n-1}):= (-1)^n2 (n-\tfrac{1}{2}) \cdots
\tfrac{1}{2} \,\text{Tr}\,(F [F, a_0][F,
a_1]\cdots [F, a_{2n-1}]).$$

\begin{definition}  The   Connes-Chern character of  an odd p-summable
Fredholm module $(\mathcal{H}, F)$ over an algebra $\A$ is the class of
the  cyclic cocycle  {\rm $\text{Ch}^{2n-1}$}
in the odd periodic cyclic cohomology group $HP^{odd}(\A)$.
\end{definition}
By the above Proposition, the class of $\text{Ch}^{2n-1}$ in $HP^{odd}(\A)$
is independent of the choice of $n$.\\

  Let us compute the character of the
Fredholm module of  the above Example  with $\A= C^{\infty} (S^1)$.
By the above definition, $\text{Ch}^1\, (\mathcal{H},\, F) = [\varphi_1]$ is
the class of the following cyclic 1-cocycle in $HP^{odd} (\A):$
$$ \varphi_1 \, (f_0, f_1)= \text{Tr} \, (F [F, f_0] [F, f_1]).$$
One can identify this cyclic cocycle with a  local
formula.
We claim that
$$ \varphi_1 \, (f_0, f_1) = \frac{4}{2\pi i}\, \int f_0 df_1, \quad
\text{for all} \,\, f_0, f_1 \in \A.$$
By linearity,  It suffices
to check this relation   for     basis elements $f_0 =e_m, f_1=e_n$ for all
$m, n \in \mathbb{Z}$, which is easy to do.

The  duality, that is,  the  bilinear pairing,  between $K$-theory and $K$-homology
is defined through the Fredholm index. More precisely there is  an
{\it index pairing} between odd (resp. even)  Fredholm modules over
$\A$ and the algebraic $K$-theory group  $K_1^{\text{alg}}(\A)$ (resp. $K_0 (\A)$). We shall
describe it only in the odd case at hand. Let
$(\mathcal{H}, \, F)$ be an odd Fredholm module over $\A$ and let $U\in
\A^{\times}$ be an invertible element in $\A$. Let $P=
\frac{F+1}{2}: \mathcal{H} \to \mathcal{H}$ be the projection operator defined by $F$.
One checks that the operator
$$ PUP: P \mathcal{H} \to P \mathcal{H}$$
 is a Fredholm operator. In fact, using the compactness of commutators $[F, a]$, one checks that
 $PU^{-1}P$ is an inverse for $PUP$
modulo compact operators, which of course implies that $PUP$ is a Fredholm operator.
 The index pairing  is then defined as
$$ \langle (\mathcal{H}, \,F),\, [U] \rangle: =\text{index}\,  (PUP),$$
where the index on the right hand side is the Fredholm index. If
the invertible $U$ happens to be in $M_n(\A)$ we can apply this
definition to the algebra $M_n(\A)$ by noticing that
$(\mathcal{H}\otimes \mathbb{C}^n, F\otimes 1)$ is a Fredholm
module over $M_n(\A)$ in a natural way. The resulting map can be
shown to induce a well defined additive map
$$ \langle (\mathcal{H}, \,F),\, - \rangle: K_1^{\text{alg}}(\A) \to \mathbb{C}.$$
Notice that this map is purely topological in the sense that
to define it     we did not have to impose any
finite summability, i.e.,  smoothness,  condition on the Fredholm module.

Going back to our example and choosing
$f: S^1\to GL_1(\mathbb{C})$ a  continuous function on $S^1$ representing  an
element of $K_1^{\text{alg}} (C(S^1))$,   the index pairing
$\langle [(\mathcal{H},\, F)], \, [f] \rangle = \text{index} (PfP)$ can be explicitly calculated.
 In fact in this case a simple homotopy
argument gives  the index of the {\it Toeplitz
operator}  $PfP: P\mathcal{H}  \to P\mathcal{H}$
in terms of the   winding number of $f$ around the origin:
$$ \langle [(H,\, F)], \, [f] \rangle =  -W (f, 0).$$
Of course, when  $f$ is smooth  the winding
number can be computed by integrating the 1-form $ \frac{1}{2\pi
i}\frac{dz}{z}$ over the curve  defined by  $f$:
$$W (f, 0)=\frac{1}{2\pi i}\int f^{-1} df =\frac{1}{2\pi i} \varphi (f^{-1}, f)$$
where $\varphi$ is the cyclic 1-cocycle on $C^{\infty} (S^1)$ defined by
$\varphi (f, g)=\int fdg$. This is a special case of a very
general index formula  proved by  Connes \cite{ac85} in a  fully noncommutative situation:

\begin{prop}\label{cindexfodd} Let $(\mathcal{H},\, F)$ be an odd  $p$-summable Fredholm
module over an algebra\ $\A$ and let $2n-1$ be an odd integer such that $2n\geq p$. If $u$ is an invertible
element  in $\
A$ then
\rm{$$ \text{index} \, (PuP)= \frac{(-1)^n}{2^{2n}} \varphi_{2n-1}
(u^{-1}, u, \dots, u^{-1}, u),$$}
where the cyclic cocycle $\varphi_{2n-1}$ is defined by
\rm{$$ \varphi_{2n-1} \,(a_0, a_1, \dots, a_{2n-1})= \text{Tr}\,  (F [F, a_0][F,
a_1]\cdots [F, a_{2n-1}]).$$}
\end{prop}

The above index formula can be expressed in a   more conceptual
manner once Connes' Chern character in $K$-theory is introduced.
In \cite{ac80, ac85}, Connes shows that the
 Chern-Weil definition of Chern character on topological  $K$-theory
 admits a vast generalization to a noncommutative
 setting. For a  noncommutative
algebra $\A$
and  each integer $n\geq 0,$ he defined  pairings  between cyclic cohomology and $K$-theory:
\begin{equation} \label{ccp}
 HC^{2n}(\A) \otimes K_0 (\A) \longrightarrow  \mathbb{C},\quad \quad
HC^{2n+1}(\A) \otimes K_1^{\text{alg}}(\A)  \longrightarrow  \mathbb{C}.
\end{equation}
These pairings are  compatible with the periodicity
operator $S$ in cyclic cohomology in the sense that
$$\langle [\varphi], \, [e] \rangle = \langle  S[\varphi], [e] \rangle,$$
for all cyclic cocycles $\varphi$ and $K$-theory classes $[e]$,
and thus induce a pairing
$$HP^{i}(\A) \otimes K_i^{\text{alg}}(\A) \longrightarrow \mathbb{C}, \quad i=0, 1$$
between periodic cyclic cohomology and $K$-theory.

We briefly recall its definition.  Let $\varphi$ be a cyclic $2n$-cocycle  on $\A$ and
let  $e \in
M_k(\A)$ be an idempotent representing a class in $K_0(\A)$. The pairing
$\label{paircyc0} HC^{2n} (\A) \otimes K_0 (\A) \longrightarrow
\mathbb{C}$
is defined by 
\begin{equation}\label{paircyc} \langle [\varphi], \, [e] \rangle = (n!)^{-1}\, \tilde{\varphi}
(e, \dots,e),
\end{equation}
where $\tilde{\varphi}$ is the `extension' of $\varphi$ to $M_k(\A)$
defined by 
the formula
\begin{equation}\label{extcoc}
\tilde{\varphi} (m_0\otimes a_0, \dots, m_{2n}\otimes a_{2n})=\text{tr}(m_0\cdots
m_{2n})\varphi (a_0, \dots, a_{2n}).
\end{equation}
  It can be shown that $\tilde{\varphi}$ is  a cyclic $n$-cocycle as well.

The formulas in the {\it odd case} are as follows. Given an invertible
matrix $u\in M_k(\A)$, representing a class in $K^{\text{alg}}_1 (\A)$,   and an odd
cyclic $(2n-1)$-cocycle $\varphi$ on $\A$, the pairing is given by
\begin{equation}\label{oddpair}
\langle [\varphi], \, [u]\rangle:  = \frac{2^{-(2n+1)}}{(n-\frac{1}{2})\cdots \frac{1}{2}}\,\tilde{\varphi}(u^{-1}-1, u-1,
\dots, u^{-1}-1, u-1).
\end{equation}
Any cyclic cocycle  can be represented by a
{\it normalized} cocycle for which  $\varphi (a_0, \dots, a_n)=0$ if
$a_i=1$ for some $i$. When $\varphi $ is normalized, formula
\eqref{oddpair} reduces to a particularly simple form:
\begin{equation}\label{odnopair}
\langle [\varphi], \, [u]\rangle =
\frac{2^{-(2n+1)}}{(n-\frac{1}{2})\cdots \frac{1}{2}}\,\tilde{\varphi}(u^{-1}, u, \dots, u^{-1}, u).
\end{equation}

Using the    pairing $ HC^{2n-1} (\A) \otimes K^{\text{alg}}_1 (\A) \to
\mathbb{C} $  and the definition of $\text{Ch}^{2n-1} \, (H, \, F),$
the above index formula 
in Proposition \eqref{cindexfodd}  can be written as
\begin{equation}\label{indf1} \text{index} \, (PuP)=
\langle \text{Ch}^{2n-1} \, (H, \, F), \, [u]\rangle,
\end{equation}
or in its stable form
$$ \text{index} \, (PuP)=  \langle \text{Ch}^{odd} \, (H, \, F), \, [u]\rangle.$$
 This equality   amounts to the equality
   {\it Topological Index = Analytic Index} in a fully noncommutative
   setting.

An immediate consequence of the index formula  \eqref{indf1} is an integrality theorem for
numbers defined by the right hand side of  \eqref{indf1}.  This should be compared
with classical integrality
results for topological invariants of manifolds that are established through the Atiyah-Singer
index theorem. An early nice application was Connes' proof of the idempotent conjecture for group
$C^*$-algebras of {\it free groups} in \cite{ac85}. Among other applications
I should mention a mathematical treatment of integral quantum Hall effect, and most recently
 to quantum
computing in the work of Mike Freedman and collaborators \cite{freed}.

\section{Cyclic modules}

 \index{cyclic module} With the introduction of the  {\it cyclic
 category}
 $\Lambda$ in \cite{ac83b}, Connes took another major step in conceptualizing
 and generalizing  cyclic cohomology  far beyond its
 original inception. We already saw in the last section three  different
 definitions
  of the cyclic cohomology of an algebra through explicit complexes.  The original motivation of \cite{ac83b} was
 to define  the cyclic cohomology of algebras as a derived functor. Since the category
 of algebras and algebra homomorphisms is not  an additive category,
 the standard
abelian homological algebra is not applicable here. Let $k$ be a unital
commutative ring. In \cite{ac83b},  an abelian  category $\Lambda_k$  of
 cyclic $k$-modules is defined that can be thought of as the   `abelianization' of the category of
$k$-algebras.
  Cyclic cohomology is then shown to be the derived
functor of the  {\it functor of traces},  as we shall explain in this section.
More generally Connes defined     the notion of  a {\it
cyclic
 object} in an abelian category and  its cyclic cohomology \cite{ac83b}.

 Later developments
 proved that this extension of cyclic cohomology was of great significance. Apart from earlier applications,  we should mention the
 recent
 work \cite{cocoma}  where the abelian category of
 cyclic modules plays a role similar to that of
 the  category of motives for noncommutative geometry.
   Another  recent example is the cyclic
 cohomology of  Hopf algebras \cite{como2, como3, hkrs1, hkrs2},  which cannot be defined as the cyclic cohomology of
an algebra or a coalgebra but only   as the cyclic cohomology of a cyclic
module naturally attached to the given Hopf algebra and a coefficient
system (see the last  section for more on Hopf cyclic cohomology). Let us briefly sketch the definition of the cyclic category $\Lambda$.

Recall that the     simplicial category $\Delta$ is a  small  category whose
objects are the totally ordered sets
 $$[n]=\{0< 1< \cdots < n\}, \quad \quad n=0, 1, 2, \dots,$$
 and whose   morphisms $f: [n] \to [m]$ are  order preserving, i.e.
 monotone non-decreasing, maps
 $f: \{0, 1, \dots, n\}\to \{0, 1, \dots, m \}.$ Of particular
 interest among the morphisms of $\Delta$ are  {\it faces} $\delta_i$
 and {\it degeneracies} $\sigma_j$,
 $$ \delta_i: [n-1]\to [n], \quad i=0, 1, \dots, n,$$
$$\sigma_j: [n+1]\to [n], \quad
 \quad j= 0, 1,  \dots, n. $$
 By definition $\delta_i$ is the unique injective morphism missing $i$ and
 $\sigma_j$ is the unique surjective  morphism identifying $j$ with $j+1$.

The  {\it cyclic category}  $\Lambda$ has the same set of objects as $\Delta$
and in fact contains $\Delta$ as a subcategory.  Morphisms of   $\Lambda$ are generated by simplicial morphisms
 and new morphisms  $\tau_n : [n] \to [n], n \geq 0,$ defined by $\tau_n (i)= i+1$ for $0\leq i <n$ and $\tau_n (n)=0.$ We have
  the following
extra relations:
\begin{eqnarray*}
\tau_n\delta_{i} =\delta_{i-1} \tau_{n-1}, \quad  \tau_n \delta_0 = \delta_{n}, & &\quad  1\le i\le n,\\
 \tau_n \sigma_i = \sigma _{i-1} \tau_{n+1}, \quad \tau_n \sigma_0 =\sigma_n \tau_{n+1}^2 & & \quad 1\le i\le n,\\
\tau_n^{n+1} = \mbox{id.} \quad & &
\end{eqnarray*}
It can be shown that the classifying space $B \Lambda$ of the
small category $\Lambda$ is homotopy equivalent to the classifying
space of the circle $S^1$ \cite{ac83b}.

 A {\it cyclic object} in a category $\mathcal{C}$ is
a  functor $\Lambda^{\text{op}}  \rightarrow \mathcal{C}$. A {\it cocyclic object} in
 $\mathcal{C}$ is a  functor $\Lambda \rightarrow \mathcal{C}$.
For any commutative unital  ring $k$, we denote the category of cyclic
  $k$-modules by $\Lambda_k $. A morphism of cyclic $k$-modules is a natural transformation between the corresponding
  functors.  It is clear that
  $\Lambda_k$ is an abelian category.  More generally,  if $\mathcal{A}$ is an abelian category then the category
  $\Lambda \mathcal{A}$  of cyclic objects in $\mathcal{A}$ is itself
  an abelian category.

Let $\text{Alg}_k$ denote the category  of  unital $k$-algebras and  unital
algebra homomorphisms. There is a functor
$$ \natural : \text{Alg}_k \longrightarrow \Lambda_k, \quad A\mapsto A^ \natural, $$
defined by
$$A_n^\natural =A^{\otimes(n+1)},  \quad \quad n\geq 0,$$
  with  face,  degeneracy
and cyclic operators  given  by
\begin{eqnarray*}
\delta_i(a_0 \otimes a_1\otimes \dots \otimes a_n)&=&a_0 \otimes \dots \otimes
a_{i}a_{i+1}\otimes \dots \otimes a_n,\\
\delta_n(a_0 \otimes a_1\otimes\dots \otimes a_n)&=&a_na_0 \otimes a_1 \otimes
\dots  \otimes a_{n-1},\\
\sigma_i(a_0 \otimes a_1\otimes\dots \otimes a_n)&=&a_0 \otimes\dots \otimes
a_i \otimes  1 \otimes \dots\otimes
 a_n,\\
\tau_n(a_0 \otimes a_1\otimes\dots \otimes a_n)&=& a_n \otimes a_0\otimes \dots \otimes a_{n-1}.
\end{eqnarray*}
A unital algebra map $f: A \to B$ induces a morphism of cyclic modules
$f^{\natural}: A^ \natural \to B^ \natural$ by $f^\natural (a_0\otimes
\cdots \otimes a_n)= (f(a_0)\otimes
\cdots \otimes f(a_n)).$

This functor $\natural$ embeds the non-additive category of $k$-algebras into the abelian
category of cyclic $k$-modules. A first main observation of
\cite{ac83b} is that
$$\text{Hom}_{\Lambda_k} \,(A^ \natural, \, \, k^ \natural) \simeq T(A),$$
where $T(A)$ is the space of traces from $A \to k$. To a trace $\tau$ one associate  the cyclic map $(f_n)_{n\geq 0}$,  where
$$ f_n(a_0\otimes a_1\otimes \cdots \otimes a_n)=\tau (a_0a_1\cdots
a_n), \quad \, n \geq 0.$$ It can be easily shown that this
defines a one to one correspondence.

Now we can state the following fundamental theorem of  Connes \cite{ac83b} which
greatly extends the above observation and shows that cyclic cohomology is  a derived functor, in fact an Ext functor,
provided that we  work in the category of cyclic modules:
\begin{theorem} \label{cycderfunc}Let $k$ be a field of characteristic zero. For any unital  $k$-algebra $A$, there is a    canonical
isomorphism
\rm{$$ HC^n(A) \simeq \text{Ext}^n_{\Lambda_k} (A^ \natural, \, k^ \natural),
\quad \quad \text{for all} \, \, n\geq 0.$$}
\end{theorem}

Apart from their applications in the study of cyclic cohomology of
algebras and Hopf algebras (about the latter see the next section), cyclic
modules have also come to play an important role in applications of
noncommutative geometry to number theory. They play a role similar to that  of motives in
algebraic geometry. Let me  briefly explain this point.

The program outlined by Connes, Consani and Marcolli in \cite{cocoma} aims at creating an environment
where something like  Weil's    proof of the Riemann hypothesis for function fields
can be repeated in the
characteristic zero case. Among  other things, they  produce an analogue of
the  Frobenius automorphism  in characteristic zero in this paper.    Since
Connes' trace formula 
    is over the
 noncommutative  {\it ad\`{e}les class space} \cite{ac99}, the geometric setting is that of noncommutative
 geometry and they must go far beyond
  what is done so far in noncommutative geometry  and import many ideas
 from modern algebraic geometry to noncommutative geometry.
To achieve this, as a first step,  good analogues of {\it \'{e}tale
cohomology}, the
 {\it category of motives}, and {\it correspondences} in
noncommutative geometry must be introduced. Happily it turns out that
  Connes' category of
 cyclic modules and the closely related bivariant cyclic homology, as well as $KK$-theory, are quite
 useful in this regard.

 The construction of the Frobenius in characteristic zero
follows a very general process that combines  cyclic homology with
quantum statistical mechanics in a novel way. Starting from a pair $(A, \varphi)$ of an algebra and a
state $\varphi$
(a noncommutative space endowed with a `probability measure'), they proceed by  invoking the
canonical one-parameter group of automorphisms $\sigma =\sigma_{\varphi}$ and consider the extremal equilibrium
states
$\Sigma_{\beta}$ at inverse temperatures $\beta >1$. Under suitable conditions,  there is
an algebra map $$\rho : A \rtimes_{\sigma} \mathbb{R} \to
\mathcal{S}(\Sigma_{\beta}\times \mathbb{R}^*_+)\otimes \mathcal{L},$$
where $\mathcal{L}$ denotes the algebra of trace class operators. The cyclic module $D(A, \varphi)$ is defined as
the cokernel of the induced map by $\text{Tr}\circ \rho$ on the
cyclic modules of these two algebras. The dual multiplicative group $\mathbb{R}^*_+$
acts on $D(A, \varphi)$ and, in examples coming from number theory,
replaces  Frobenius in characteristic zero. The three steps involved in the
construction of $D(A, \varphi)$ are called {\it cooling}, {\it distillation},
and {\it dual action} in the paper.

 A remarkable property of the  cyclic category $\Lambda$, not shared  by the simplicial category,
  is its {\it self-duality} in the sense that there is a
natural isomorphism  of categories $\Lambda \simeq \Lambda^{\text{op}}$ \cite{ac83b}.   Roughly  speaking, the   duality
functor $\Lambda^{\text{op}} \longrightarrow
\Lambda $ acts as the  identity on objects
of $\Lambda$ and exchanges face and degeneracy operators while sending
the cyclic operator to its inverse.
Thus to a cyclic (resp. cocyclic) module one can
associate a cocyclic (resp. cyclic) module by applying the duality
isomorphism. This duality plays an important role in Hopf cyclic
cohomology.

\section{The local index formula and beyond}

 In practice, computing   Connes-Chern characters
 defined by formulas like \eqref{oddcc}
 is rather difficult since they  involve the
ordinary operator trace and are  non-local.  Thus one needs to compute
the class of this cyclic cocycle
by a {\it local formula}. This is rather similar to passing  from
the   McKean-Singer formula  for the index of an elliptic operator to a local
cohomological formula involving integrating  a  locally defined differential form, i.e., the
Atiyah-Singer index formula.  The solution of this problem was arrived
at in two stages. First,
in \cite{acbook}, Connes gave a partial answer by giving  a local  formula for the {\it Hochschild
class} of the Chern
character, and then Connes and
Moscovici gave a formula that captures  the full cyclic cohomology class of the character
by a local formula \cite{como1}. Broadly speaking,  the ideas involved
amount to going from {\it noncommutative  differential topology} to {\it noncommutative
 spectral geometry},  and need the introduction of  two new concepts.

In the first place,   a  noncommutative analogue of integration was  found by  Connes by replacing the operator
trace by
the Dixmier trace \cite{acaction}, and, secondly,   one  refines  the topological notion
of Fredholm module by the metric notion of a {\it spectral triple}, or $K$-cycles as
they were originally named in \cite{acbook}.
Developing the necessary tools to handle this local index formula,  shaped,
more or less,  the second stage of the development of noncommutative
geometry after the appearance of the landmark papers \cite{ac85}. One can say that while in its first stage
noncommutative geometry was influenced
by differential and algebraic topology,   especially  index theory, the Novikov conjecture
and the Baum-Connes conjecture,  in this second stage it was  chiefly informed by spectral geometry.

We start with a quick
review of the {\it Dixmier trace} and the {\it noncommutative integral}, following  \cite{acbook} closely.
For a compact  operator $T$,  let $\mu_n(T), n=1, 2, \dots,$ denote the sequence of eigenvalues
of $|T|=(T^*T)^{\frac{1}{2}}$ written in  decreasing order. Thus, by the  minimax principle, $\mu_1 (T) = ||T||$, and in general
$$ \mu_n (T) = \text{inf} \, \,||T|_V||, \quad n\geq1,$$
where the infimum is over the set of subspaces of codimension $n-1$, and
$T|_V$ denotes  the restriction of $T$ to the subspace $V$.
The natural domain of the Dixmier trace is the set of  operators
$$\mathcal{L}^{1,\infty } (\mathcal{H}):= \{T \in
\mathcal{K}(\mathcal{H}); \, \sum_1^N \mu_n (T)=O \,(\text{log} N )\}.$$
Notice that trace class operators are automatically in $\mathcal{L}^{1,\infty} (\mathcal{H})$.
The Dixmier trace of an operator $T \in \mathcal{L}^{1,\infty}
(\mathcal{H})$
  measures the {\it logarithmic divergence} of
its ordinary trace. More precisely,  we are interested in  taking some kind
of limit of the bounded sequence
$$\sigma_N (T)= \frac{\sum_1^N \mu_n(T)}{\text{log} N}$$
as $N\to \infty$.
The problem of course is that, while by our assumption the sequence is bounded,  the usual limit may
not exists and must be replaced by a carefully chosen `generalized
limit'.

To this end, let $\text{Trace}_{\Lambda} (T), \Lambda \in [1, \infty)$,  be
the piecewise affine interpolation of the partial trace function
$\text{Trace}_{N} (T) =\sum_1^N \mu_n (T)$.  Recall that a state on a  $C^*$-algebra
is a non-zero positive linear functional on the algebra. Let
$\omega : C_b[e, \infty) \to \mathbb{C}$ be a normalized state on the algebra of bounded continuous functions
on  $[e, \infty)$ such that $\omega (f)=0$ for all f vanishing at $\infty$. Now, using $\omega$,  the
Dixmier trace of a positive operator $T \in \mathcal{L}^{1,\infty}
(\mathcal{H})$ is defined as
$$ \text{Tr}_{\omega}(T): = \omega (\tau_{\Lambda} (T)),$$
where
$$\tau_{\Lambda} (T)=  \frac{1}{\text{log} \, \Lambda} \int_e^{\Lambda} \frac{ \text{Trace}_{r} (T)}{\text{log}\,r}\frac{dr}{r}$$
 is the  Ces\`{a}ro mean  of the function $\frac{ \text{Trace}_{r} (T)}{\text{log}r}$
 over the multiplicative group $\mathbb{R}^*_+$. One then extends
 $\text{Tr}_{\omega}$ to all of $\mathcal{L}^{1,\infty}
(\mathcal{H})$ by linearity.

The resulting linear functional  $\text{Tr}_{\omega}$ is a
  positive trace on $\mathcal{L}^{(1, \infty)} (\mathcal{H})$.  It is easy to see from its definition
  that if  $T$ actually
  happens to be a trace class operator then $\text{Tr}_{\omega}(T)=0$
  for all $\omega$, i.e., the Dixmier trace is invariant under 
  perturbations by trace class operators. This is a very useful property and  makes $\text{Tr}_{\omega}$ a  flexible
  tool in computations.
   The Dixmier trace,   $\text{Tr}_{\omega}$,  in general depends on the
limiting procedure $\omega$;    however,  for the class of operators
$T$  for which   $\text{Lim}_{\Lambda \to \infty}\,  \tau_{\Lambda}(T) $
exit, it  is  independent of the choice of $\omega$ and is equal to
$\text{Lim}_{\Lambda \to \infty} \tau_{\Lambda}(T)$. One of the  main results proved in
\cite{acaction} is that if
$M$ is a closed  $n$-dimensional manifold,  $E$ is a smooth vector bundle on  $M$, $P$
is a  pseudodifferential
operator of order $-n$  acting between $L^2$-sections of $E$, and $H = L^2 (M, E),$ then
$P\in \mathcal{L}^{(1, \infty)} (\mathcal{H})$   and, for any choice of
$\omega$,
$\text{Tr}_{\omega} (P)= n^{-1}
\text{Res} (P)$. Here Res denotes Wodzicki's  noncommutative residue.
For example, if $D$ is an elliptic first order differential
operator,  $|D|^{-n}$ is a pseudodifferential operator of order $-n$ and,
for any bounded operator $a$, the Dixmier trace
$\text{Tr}_{\omega} (a |D|^{-n})$ is independent of the choice of $\omega$.

The second ingredient of the local index formula is the notion of
{\it spectral triple}  \cite{acbook}. Spectral triples provide a refinement of Fredholm modules. Going from Fredholm modules to spectral triples is
similar to going from the conformal class of a Riemannian metric to the
metric itself. Spectral triples simultaneously  provide a notion of
{\it Dirac operator} in noncommutative geometry, as well as a Riemannian
type {\it distance function} for noncommutative spaces.

To motivate the definition of a spectral triple, we recall that the Dirac operator ${D \hspace{-7pt} \slash}$ on a
compact Riemannian $\text{Spin}^c$  manifold  acts  as an unbounded
selfadjoint operator on  the Hilbert space $L^2(M,S)$ of $L^2$-spinors on the manifold $M$. If we  let
$C^{\infty}(M)$ act  on $L^2(M,S)$ by multiplication operators, then  one can check that for any smooth
function $f$, the commutator $[D,f]=Df-fD$ extends to a bounded
operator on  $L^2(M,S)$. Now  the geodesic distance $d$ on $M$ can be recovered from the following beautiful
{\it distance formula} of Connes \cite{acbook}:
\[d(p,q) = \textnormal{Sup} \{ |f(p)-f(q)|; \parallel [D,f] \parallel \leq 1\}, \quad \,\,\, \forall p,q \in M.
\nonumber \]
The triple $(C^{\infty}(M), L^2(M,S), {D \hspace{-7pt} \slash})$ is a commutative example of
a spectral triple. Its general definition, in the odd case, is as
follows.
This definition  should be compared with
Definition \eqref{fm1}.
\begin{definition}
Let $\A$ be a unital algebra. An odd spectral triple on $\A$ is a triple $(\A,
\mathcal{H}, D)$ consisting  of a Hilbert space $\mathcal{H}$, a selfadjoint unbounded operator
{\rm $D:  \text{Dom} (D) \subset \mathcal{H} \to \mathcal{H}$} with
compact resolvent, i.e.,
$(D+\lambda)^{-1} \in \mathcal{K}(\mathcal{H}),  \text{for all} \,
\lambda \notin \mathbb{R},$
 and a
representation  $\pi: \A \to \mathcal{L}(\mathcal{H})$  of $\A$ such that for
all $a\in \A$, the commutator $[D, \pi (a)]$ is defined on {\rm $ \text{Dom} (D)$} and extends to  a bounded operator on $\mathcal{H}$.
\end{definition}

The finite summability assumption  \eqref{fs11} for  Fredholm modules has a finer analogue
for spectral triples. For simplicity we shall assume  that $D$ is
invertible (in general, since $\text{Ker} \,D$ is finite dimensional, by
restricting to its orthogonal complement we can always reduce to this case).
  A  spectral
triple is called {\it finitely summable}  if for some $n\geq 1$
\begin{equation}\label{stfs} |D|^{-n} \in \mathcal{L}^{1, \, \infty}(\mathcal{H}).
\end{equation}

A simple example of an odd spectral triple is
$(C^{\infty} (S^1), L^2 (S^1), D)$, where $D$ is the unique selfadjoint
extension of the  operator
$-i \frac{d}{dx}$. Eigenvalues of $|D|$ are $|n|,  n\in \mathbb{Z}$, which shows
that, if we restrict $D$ to the orthogonal complement of constant
functions, then $|D|^{-1} \in \mathcal{L}^{1, \, \infty}(L^2(S^1)).$

Given a spectral triple $(\A, \mathcal{H}, D)$, one obtains a
Fredholm module $(\A, \mathcal{H}, F)$ by choosing $F =\text{Sign}\,
(D)= D|D|^{-1}$. Connes' Hochschild character formula gives a local
expression for the Hochschild class of the Connes-Chern character of
$(\A, \mathcal{H}, F)$ in terms of $D$ itself. For this one has to assume that the spectral  triple $(\A, \mathcal{H}, D)$ is
 {\it regular}  in the sense that  for all $a\in \A$,
$$ a \, \, \text{and} \,\, [D, a] \, \in  \cap \,\text{Dom} (\delta^k)$$
where
 the derivation $\delta$ is  given by $\delta (x)= [|D|, x].$

 Now,  assuming \eqref{stfs} holds, Connes
defines an $(n+1)-$linear functional $\varphi$ on $\A$ by
$$\varphi (a^0, a^1, \dots, a^n)= \text{Tr}_{\omega} (a^0[D, a^1]\cdots
[D, a^n]|D|^{-n}). $$
It can be shown that $\varphi$ is a Hochschild $n$-cocycle on $\A$.
We recall that a Hochschild $n$-cycle $c \in Z_n (A, A)$ is an element
$c  =\sum a^0\otimes a^1 \otimes \dots \otimes a^n \in A^{\otimes (n+1)}$ such
that its Hochschild boundary $b (c)=0$, where $b$ is defined by \eqref{hhomb1}.
The
 following result, known as {\it Connes' Hochschild character formula}, computes the Hochschild
 class of the Chern charcater by a local formula, i.e., in terms of
 $\varphi$:
\begin{theorem} Let  $(\A, \mathcal{H}, D)$ be a regular spectral
triple.
Let  {\rm $ F= \text{Sign}\, (D)$} denote the sign of $D$ and $\tau_n \in HC^n(\A)$  denote the  Connes-Chern charcater of
$(\mathcal{H}, F)$. For every $n$-dimensional Hochschild cycle $c =\sum a^0\otimes a^1 \otimes \dots \otimes a^n \in Z_n
(\A, \A)$, one has
$$\langle \tau_n, c\rangle =  \sum \varphi (a^0, a^1, \dots, a^n).$$
\end{theorem}

Identifying  the full cyclic cohomology class of the Connes-Chern character of $(\A, \mathcal{H}, D)$   by a
local formula is the content
of Connes-Moscovici's local index formula.
For this we have to assume the spectral triple satisfies  another
technical condition.
Let $\mathcal{B}$ denote the subalgebra of
$\mathcal{L}(\mathcal{H})$ generated by operators $\delta^k (a)$ and
$\delta^k ([D, a]), \, k\geq 1.$ A spectral triple is said to have a discrete {\it dimension spectrum} $\Sigma$
if $\Sigma \subset \mathbb{C}$ is discrete and for any $b\in
\mathcal{B}$ the function
$$ \zeta_b (z)= \text{Trace}(b |D|^{-z}), \, \,  \quad \text{Re}\, z > n,$$
extends to a holomorphic function on $\mathbb{C}\setminus  \Sigma$. It is further
assumed that $\Sigma$ is simple in the sense that  $\zeta_b (z)$ has only
simple poles in $\Sigma$.

The local index formula of Connes and Moscovici \cite{como1} is given by the
following Theorem (we have used the formulation in \cite{ac00}):
\begin{theorem}
1. The equality
{\rm $$ \int\!\!\!\!\!\!-P= \text{Res}_{z=0} \, \text{Trace}
(P|D|^{-z})$$}
defines a trace on the algebra generated by $\A$, $[D, \A]$, and $|D|^z,
\, z \in \mathbb{C}.$\\
2. There are  only a finite number of non-zero terms in the following
formula which defines the odd components $(\varphi_n)_{n=1, 3, \dots}$ of an odd cyclic cocycle in the $(b, B)$ bicomplex of
$\A$:
For each odd integer $n$ let
$$ \varphi_n (a^0, \dots, a^n):= \sum_{k}c_{n, k} \int\!\!\!\!\!\!- a^0[D, a^1]^{(k_1)}\cdots
[D, a^n]^{(k_n)}|D|^{-n-2|k|}$$
where $T^{(k)}:=  \nabla^k$ and $\nabla (T)= D^2T-TD^2$, $k$ is a
multi-index,
$|k|=k_1+\dots + k_n$ and
$$ c_{n, k}:= (-1)^{|k|}\sqrt{2 i} (k_1!\cdots k_n!)^{-1}((k_1+1)\cdots
(k_1+k_2+\cdots k_n))^{-1} \Gamma (|k|+\frac{n}{2}).$$
3. The pairing of the cyclic cohomology class  $(\varphi_n) \in HC^{*}(\A)$ with $K_1 (\A)$ gives the Fredholm index
of $D$ with coefficients in $K_1 (\A)$.
 \end{theorem}

As is indicated  in  part 1) of the above Theorem, a regular spectral triple necessarily defines  a trace on its underlying algebra by the formula
$a\in \A \mapsto  \int\!\!\!\!\!\!- a = \text{Res}_{z=0} \, \text{Trace} (a|D|^{-z})$. Thus,  to deal with `type III algebras'
   which carry  no non-trivial traces,  the notion of spectral triple
   must be modified.  In \cite{como7} Connes and Moscovici define a
   notion of
{\it twisted spectral triple}, where the twist is afforded by an algebra automorphism (related to the modular automorphism group). More precisely,
one postulates that there exists an
automorphism $\sigma$ of $\mathcal{A}$ such that the twisted commutators
\[[D,a]_{\sigma} := Da-\sigma(a)D \nonumber \]
are bounded operators for all $a \in \mathcal{A}$.
 They show that,  in the twisted case, the Dixmier trace
induces a twisted trace on the algebra $\mathcal{A}$, but surprisingly, under some regularity conditions,
the Connes-Chern character of the phase space lands in ordinary cyclic cohomology. Thus its pairing
with ordinary K-theory makes sense, and it can be recovered as the index of Fredholm operators.
This suggests  the significance of developing
a local index formula for twisted spectral triples, $i.e.$ finding a formula for a cocycle, cohomologous to
the Connes-Chern character in the $(b,B)$-bicomplex, which is given in terms of twisted commutators and
residue functionals. I beleive that this new theme of twisted spectral
triples, and type III noncommutative geometry in general,
will  dominate
the subject in near future.

For example, very recently a local index formula has been proved for a class of twisted spectral
triples by Henri Moscovici \cite{mosc1} that can be found in the present volume.
This class is obtained by twisting an ordinary spectral triple $(\mathcal{A}, \mathcal{H}, D)$ by
a subgroup $G$ of conformal similarities of the triple, $i.e.$ the set of all unitary operators $U \in \mathcal{U}
(\mathcal{H})$ such that $U \mathcal{A} U^*= \mathcal{A}$, and $UDU^*=\mu(U)D$, with $\mu(U) >0$. It is shown
that the crossed product algebra $\mathcal{A} \rtimes G$ admits an automorphism $\sigma$, given by the formula
$\sigma(aU)=\mu(U)^{-1}aU$, for all $a\in \mathcal{A}, U\in G$,  and $(\mathcal{A} \rtimes G, \mathcal{H}, D)$
is a twisted spectral triple. The analogue of the noncommutative residue
on the circle, 
for algebras of formal twisted pseudodifferential symbols, is constructed in \cite{fakh1}.

A very recent development related to (twisted) spectral triples is    the noncommutative Gauss-Bonnet
theorem of Connes and Tretkoff
for the  noncommutative two-torus $A_{\theta}$ \cite{contre}.
In classical geometry  a {\it spectral zeta function} is associated to the
Laplacian $\Delta _{g} = d^*d$ of a Riemann
surface with metric $g$:
\[ \zeta (s) = \sum_j \lambda_{j}^{-s}, \,\,\, \textnormal{Re}(s) > 1, \nonumber \]
where the $\lambda_j$' s are the nonzero eigenvalues of $\Delta _{g}$.  This zeta function has a meromorphic continuation with no pole at 0, and
the Gauss-Bonnet theorem for surfaces can be  expressed as
\[ \zeta(0) + \textnormal{Card} \{j | \lambda_j = 0\} = \frac{1}{12 \pi} \int_{\Sigma} R = \frac{1}{6}
\chi (\Sigma), \nonumber \]
where $R$ is the  curvature and $\chi (\Sigma)$ is the Euler-Poincar\'e characteristic.

It is this formulation of the Gauss-Bonnet theorem in
spectral terms that admits a generalization to noncommutative geometry.
Let $A_{\theta}$ denote the $C^*$-algebra of the noncommutative torus
with parameter $\theta \in \mathbb{R}\setminus \mathbb{Q}$ and let $\tau: A_{\theta}\to
\mathbb{C}$ denote its faithful  normalized trace.
One can define an inner product
\[ \langle a, b \rangle = \tau (b^*a), \,\,\, a,b \in A_{\theta}, \nonumber \]
and complete $A_{\theta}$ with respect to this inner product to obtain a Hilbert space $\mathcal{H}_0$.
More generally, for any smooth selfadjoint element $h=h^* \in
A_{\theta}$ one defines an inner product $\langle a, b \rangle_{\varphi} = \tau
(b^*a e^{-h})$, where the  positive linear functional $\varphi = \varphi_h$ is defined by
$$\varphi(a)= \tau (ae^{-h}), \,  \quad a\in A_{\theta}.$$
  Let $\mathcal{H}_{\varphi}$ denote the completion of $A_{\theta}$ with respect to this
conformally equivalent metric.

Using the canonical derivations $\delta_1$ and $\delta_2$ of  $A_{\theta}$, one can introduce
a complex structure on $A_{\theta}$ by defining
\[ \partial = \delta_1 + i \delta_2, \,\,\, \partial^*=  \delta_1 - i \delta_2. \nonumber \]
These operators can be considered as unbounded operators on $\mathcal{H}_0$ and $\partial^*$ is
the adjoint of $\partial$. Then the {\it unperturbed Laplacian} on $A_{\theta}$ is given by
\[\Delta= \partial^* \partial = \delta_1^2 + \delta_2^2. \nonumber \]
In general we can consider the unbounded operator  $\partial = \delta_1
+ i \delta_2:\mathcal{H}_{\varphi} \to \mathcal{H}^{(1, 0)}$, where $\mathcal{H}^{(1, 0)}$ is the
completion of the linear span of elements of the form $ a\partial b$ with $a, b \in  A_{\theta}^{\infty}$. Let $\partial^*_{\varphi}$ denote
its adjoint. Then the Laplacian for   the conformally equivalent metric
$\langle a, b \rangle_{\varphi}$ is given by $\Delta'= \partial^*_{\varphi} \partial.$

In \cite{contre}, Connes  and Tretkoff
 show that the value at 0 of the zeta function associated to this Laplacian
 $\Delta'$ is an
invariant of the conformal class of the metric on $A_{\theta}$, i.e. of $h$.   A natural problem here is to extend  this result
by considering the most general complex structure on $A_{\theta}$ of the form
$\partial = \delta_1 + \tau \delta_2$, where $\tau$ is a complex number
with $\textnormal{Im}(\tau) > 0 $. This problem is now solved in full
generality in \cite{fakh2}.

\section{Hopf cyclic cohomology}
A major development in cyclic cohomology theory in the last ten years
was the
introduction of  {\it Hopf cyclic cohomology} for Hopf algebras by Connes and Moscovici
\cite{como2}.
As we saw in Section 5, the local
index formula  gives the Connes-Chern character
of a regular
spectral triple $(\A, \mathcal{H}, D)$ as a  cyclic cocycle in the $(b, B)$-bicomplex of the algebra $\A$.
For  spectral
triples of interest in transverse geometry \cite{como2}, this cocycle is {\it differentiable} in
the sense that it is in the
image of the Connes-Moscovici characteristic map $\chi_{\tau}$ defined below \eqref{cmcmap}, with
$H=\mathcal{H}_1$ a Hopf algebra and $\A=\mathcal{A}_{\Gamma}$, a noncommutative algebra,
whose definitions
 we shall
recall  in this section.  To identify
this cyclic cocycle  in terms of characteristic classes of foliations, they realized that
it would be
extremely helpful to
 show that it is the image of a polynomial in some universal cocycles for a
 cohomology theory for a universal Hopf algebra,  and this gave birth
 to Hopf cyclic cohomology and to the universal Hopf algebra $H=\mathcal{H}_1$.  This is
similar to the situation for classical characteristic classes of manifolds, which are pullbacks of group
cohomology classes.

The Connes-Moscovici characteristic map can be formulated in general terms as follows. Let $H$ be
a Hopf algebra acting as {\it quantum symmetries} of an algebra $\A$,
i.e.,
$\A$ is a left $H$-module, and the
algebra structure of $\A$ is compatible with the coalgebra structure of
$H$ in the sense that the multiplication  $\A \times \A \to \A$ and the unit map $\mathbb{C} \to \A$ of
$\A$ are  morphisms of  $H$-modules. A common terminology to describe this situation
is to say that $\A$ is a left $H$-module algebra. Using Sweedler's notation for
the coproduct of $H$, $\Delta (h)=h^{(1)} \otimes h^{(2)}$ (summation is understood), this latter
compatibility condition can be expressed as
$$ h(ab)=h^{(1)}(a)h^{(2)}(b), \quad \text{and}\quad h(1)=\varepsilon (h)1,$$
for all $h \in H$ and $a, b \in \A$.   In general one
should think of such  actions  of Hopf algebras  as the noncommutative geometry analogue  of the action of
differential operators on a manifold.

It is also important to
extend  the   notion of  trace to  allow {\it twisted traces}, such  as
 KMS states in quantum statistical mechanics, as well as the idea of {\it invariance} of a  (twisted) trace. The general setting introduced  in \cite{como2} is the
following.
 Let $\delta: H \to \mathbb{C}$ be a character of
$H$, i.e. a unital algebra map,  and $\sigma \in H$ be a grouplike element, i.e. it satisfies
 $\Delta \sigma = \sigma \otimes \sigma$.   A
linear map
$\tau :A \rightarrow \mathbb{C}$ is called $\delta$-{\it invariant} if for
all $h \in H$ and $a\in \A$,
$$ \tau (h(a))=\delta (h) \tau (a),$$
and is called a $\sigma$-{\it trace} if for all $a, b$ in $\A$,
$$ \tau (ab)=\tau (b \sigma (a)).$$
 Now for $a, b \in \A,$ let
$$\langle a, \, b \rangle:=\tau (ab).$$
 Let $\tau $ be a $\sigma$-trace on $\A$. Then $\tau$ is
$\delta$-invariant if and only if
the {\it integration by parts formula} holds. That is,  for all $h \in H$ and $a, b \in A,$
\begin{eqnarray} \label{ibpf}
\langle h(a), \, b \rangle = \langle a, \,
\widetilde{S}_{\delta}(h)(b)\rangle.
\end{eqnarray}
Here $S$ denotes the antipode of $H$ and  the
  $\delta$-{\it twisted antipode}
$\widetilde{S}_{\delta}: H \rightarrow H$ is defined by $\widetilde{S}_{\delta}=\delta *
S$, i.e.
$$\widetilde{S}_{\delta}(h)= \delta (h^{(1)})S(h^{(2)}),$$
for all $h \in H$. Loosely speaking, \eqref{ibpf} says that the formal adjoint of the differential operator $h$ is
$\widetilde{S}_{\delta}(h)$.
 Following \cite{como2, como3}, we
say that $(\delta, \sigma)$ is a {\it modular pair} if $\delta
(\sigma)=1$, and a {\it modular pair in involution} if in addition
we have
$$ \widetilde{S}_{\delta}^2 (h)=\sigma
h\sigma^{-1},$$
for all $h$ in $H$. The importance of this notion will become clear in
the next paragraph.

 Now, for each $n\geq 0$,  the  Connes-Moscovici {\it
characteristic map}
\begin{eqnarray}\label{cmcmap} \chi_{\tau}: H^{\otimes n}\longrightarrow C^n(\A),
\end{eqnarray}
 is defined by
\begin{eqnarray*}
\chi_{\tau}(h_1\otimes \cdots \otimes h_n)(a_0\otimes \cdots \otimes a_n)=
\tau (a_0 h_1(a_1)\cdots h_n(a_n)).
\end{eqnarray*}
Notice that the right hand side of \eqref{cmcmap} is the cocyclic module that (its cohomology) defines the cyclic cohomology
of the algebra $\A$. The main question about \eqref{cmcmap} is whether one  can    promote the collection of linear
spaces $\{H^{\otimes n}\}_{n\geq 0}$ to a  cocyclic module
such that the characteristic map $\chi_{\tau}$ turns into a
morphism of cocyclic modules.  We recall that the face, degeneracy,
and cyclic operators for  $\{ C^n(\A)\}_{n\geq 0}$ are
defined by:
\begin{eqnarray*}
\delta_i \varphi (a_0, \dots, a_{n+1})&=&\varphi (a_0, \dots,
a_ia_{i+1}, \dots, a_{n+1}), \quad i=0, \dots, n,\\
\delta_{n+1} \varphi (a_0, \dots, a_{n+1})&=&\varphi (a_{n+1}a_0,
a_1, \dots, a_n),\\
\sigma_i\varphi (a_0, \dots, a_{n}) &=&\varphi (a_0, \dots,a_i, 1,
\dots, a_{n}),\quad i=0, \dots, n,\\
\tau_n \varphi (a_0, \dots, a_n)& =&\varphi (a_n, a_0, \dots, a_{n-1}).
\end{eqnarray*}

 The
relation
$h(ab)=h^{(1)}(a)h^{(2)}(b)$
shows that,   in order for  $\chi_{\tau}$ to be compatible with  face
operators, the face operators $\delta_i$ on $H^{\otimes n}$, at least for $0\leq i<n,$
must involve
 the coproduct of
$H$. In fact if we define, for $ 0\leq i \leq n$,
$\delta^n_i: H^{\otimes n} \to H^{\otimes (n+1)}$, by
\begin{eqnarray*}
\delta_0(h_1\otimes
\cdots \otimes h_n)&=& 1\otimes h_1\otimes
\cdots \otimes h_n,\\
 \delta_i(h_1\otimes
\cdots \otimes  h_n)&=& h_1\otimes \cdots \otimes h_i^{(1)}\otimes h_i^{(2)}\otimes \cdots
\otimes h_n,\\
\delta_{n+1}(h_1\otimes \cdots \otimes h_n)&=&h_1\otimes
\cdots \otimes  h_n\otimes \sigma,
\end{eqnarray*}
then we have, for all $i=0, 1, \dots,n$
$$\chi_{\tau}\delta_i =\delta_i \chi_{\tau}.$$ Notice that the last
relation is a consequence of the $\sigma$-trace property of $\tau$.
Similarly, the relation $h(1_A)=\varepsilon (h) 1_A$
shows that the degeneracy operators on $H^{\otimes n}$ should
involve the counit of $H$. We thus define
$$\sigma_i(h_1\otimes \cdots \otimes h_n)=h_1 \otimes \cdots \otimes \varepsilon
(h_i)\otimes \cdots \otimes h_n.$$

It is very hard, on the other hand,  to come up with a correct  formula  for  the
 {\it cyclic operator} $\tau_n : H^{\otimes n} \to H^{\otimes n}$.
 Compatibility with $\chi_{\tau}$ demands that
$$ \tau (a_0  \tau_n(h_1\otimes \cdots \otimes h_n)(a_1\otimes \cdots \otimes a_n))=
\tau (a_n  h_1(a_0)h_2(a_1)\cdots h_n(a_{n-1})),$$
for all $a_i$'s and $h_i$'s.
For $n=1$, the  integration by parts formula \eqref{ibpf}, combined with the $\sigma$-trace property of $\tau$,
shows that
$$\tau (a_1 h(a_0))=\tau (h(a_0)\sigma (a_1))=\tau (a_0\tilde{S}_{\delta}(h)\sigma (a_1)).$$
This suggests that we should define $\tau_1: H \to H$ by
$$\tau_1(h)= \tilde{S}_{\delta}(h)\sigma.$$
Note that the required cyclicity condition for $\tau_1,$
$\tau_1^2=1,$ is equivalent to the involution  condition
$\widetilde{S}_{\delta}^2 (h)=\sigma
h\sigma^{-1}$ for the  pair $(\delta, \sigma)$.
This line of reasoning can be extended to all $n\geq 0$ and gives us:
\begin{eqnarray*}
\tau (a_n h_1(a_0) \cdots h_n(a_{n-1}))& =&\tau (
h_1(a_0) \cdots h_n(a_{n-1})\sigma (a_n))\\
&=&\tau(a_0 \tilde{S}_{\delta}(h_1)(h_2 (a_1)\cdots h_n(a_{n-1})\sigma (a_n)))\\
&=& \tau(a_0 \tilde{S}_{\delta}(h_1)\cdot(h_2 \otimes \cdots \otimes h_n\otimes \sigma)
(a_1\otimes \cdots \otimes a_n)).
\end{eqnarray*}
This suggests that the {\it Hopf-cyclic operator} $\tau_n : H^{\otimes
n}\to  H^{\otimes n}$ should be defined as
$$\tau_n (h_1\otimes \cdots \otimes h_n)=\tilde{S}_{\delta}(h_1)\cdot(h_2 \otimes \cdots \otimes
h_n\otimes \sigma),$$
where $\cdot$ denotes the diagonal action defined by
$$h\cdot (h_1\otimes \cdots \otimes h_n):= h^{(1)}h_1\otimes
h^{(2)}h_2\otimes \cdots \otimes h^{(n)}h_n.$$

 The
remarkable fact, proved by Connes and Moscovici \cite{como2, como3},
is that endowed with the above face, degeneracy, and cyclic
operators, $\{H^{\otimes n} \}_{n\geq 0}$ is a cocyclic module. The
proof is a very clever and complicated tour de force of Hopf algebra
identities. 

The resulting cyclic cohomology groups, which depend on the choice of a modular pair in involution $(\delta, \sigma)$,
 are denoted by
$HC^n_{(\delta, \sigma)}(H)$, $n=0, 1,\dots.$ The characteristic map \eqref{cmcmap} clearly induces a map between corresponding cyclic cohomology
groups
$$\chi_{\tau}: HC^n_{(\delta, \sigma)}(H) \to HC^n(\A).$$
Under this map  Hopf cyclic cocycles  are mapped  to cyclic cocycles on $\A$. Very many of the interesting cyclic cocycles in noncommutative geometry
are obtained in this fashion. Using the above discussed cocyclic module structure of $\{H^{\otimes n} \}_{n\geq 0}$, we see that
a Hopf  cyclic $n$-cocycle is
an element $x\in H^{\otimes n}$ which satisfies the relations
$$bx=0, \quad (1-\lambda)x=0,$$
where $b: H^{\otimes n} \to H^{\otimes (n+1)}$ and $\lambda:  H^{\otimes
n}\to H^{\otimes n}$ are defined by
\begin{eqnarray*}
 b(h^1\otimes \cdots \otimes h^n)&=&1\otimes h_1\otimes
\cdots \otimes h_n\\
&+&\sum_{i=1}^n (-1)^i  h_1\otimes \cdots \otimes h_i^{(1)}\otimes h_i^{(2)}\otimes \cdots
\otimes h_n\\
&+&(-1)^{n+1}h_1\otimes
\cdots \otimes  h_n\otimes \sigma, \\
\lambda (h_1\otimes \cdots \otimes h_n)&=& (-1)^n\tilde{S}_{\delta}(h_1)\cdot(h_2 \otimes \cdots \otimes
h_n\otimes \sigma).
\end{eqnarray*}

The  characteristic map \eqref{cmcmap} has its origins in Connes'
earlier work on  noncommutative differential geometry \cite{ac80}, and on
his work on the transverse fundamental class of foliations \cite{ac83c}.
In fact in these papers
 some interesting cyclic cocycles
were defined  in the context of actions of Lie algebras and
(Lie) groups. Both  examples can be shown to be special cases of the
characteristic map.
 For  example let $\A= \A_{\theta}$ denote the smooth  algebra of
 coordinates for
  the noncommutative torus with parameter $\theta \in \mathbb{R}$. 
The abelian Lie algebra $\mathbb{R}^2$ acts on $\A_{\theta}$ via canonical derivations
 $\delta_1$ and
   $\delta_2$.  
The standard trace $\tau$ on $\A_{\theta}$ is invariant under the action of $\mathbb{R}^2$, i.e.,  we
have $\tau(\delta_1(a))=\tau(\delta_2(a))=0$ for all $a \in \A_{\theta}.$
   Then one can directly check 
 that under the characteristic map \eqref{cmcmap} the two dimensional generator of the Lie algebra
 homology of $\mathbb{R}^2$ is mapped to the following   cyclic $2$-cocycle on
 $\A_{\theta}$ first defined in \cite{ac80}:
$$\varphi(a_0,a_1,a_2)=\frac{1}{2 \pi i}\tau(a_0(\delta_1(a_1)\delta_2(a_2)-\delta_2(a_1)\delta_1(a_2))).$$

  For a second example let $G$ be a discrete group and ${c}$
 be  a normalized group $n$-cocycle on $G$ with trivial coefficients. Here by normalized we mean $c(g_1, \dots, g_n)=0$ if $g_i=e$ for some $i$.
 Then one  checks that the following
 is a cyclic $n$-cocycle on the group algebra $\mathbb{C} G$~\cite{ac83c}:
 \begin{align*}
 \varphi (g_0, g_1 \dots, g_n)=\left\{ \begin{matrix} {c}(g_1, g_2\dots, g_n) \quad \text{if}~~g_0 g_1
  \dots g_n=1\\
 0  \quad \quad \quad \quad \text{otherwise}\end{matrix}\right.
 \end{align*}
After an appropriate dual version of Hopf cyclic cohomology is defined,
one can show that this cyclic cocycle can also be defined via \eqref{cmcmap}.

 The most sophisticated example of the characteristic map
 \eqref{cmcmap},  so far,  involves
 the Connes-Moscovici Hopf algebra $\mathcal{H}_1$ and its action
 on algebras of interest in transverse geometry.  In fact, as we  shall see,
 $\mathcal{H}_1$ acts as  quantum symmetries of  various objects of interest in noncommutative
 geometry, like the frame bundle of the `space' of leaves of codimension one foliations
 or the `space' of modular forms
 modulo the action of Hecke correspondences.

 To describe $\mathcal{H}_1$, let $\mathfrak{g}_{aff}$  denote  the Lie algebra of the
group of affine transformations of the line with linear basis $X$
and $Y$ and the relation $ [Y,X] = X$. Let $\mathfrak{g}$ be an
abelian Lie algebra with basis $\{\delta_n; \; n=1, 2, \dots\}$. Its
universal enveloping algebra $U(\mathfrak{g})$ should be regarded as the
continuous dual of the pro-unipotent group of orientation preserving diffeomorphisms $\varphi$ of
$\mathbb{R}$ with $\varphi (0)=0$ and $\varphi'(0)=1$.
It is easily seen that $\mathfrak{g}_{aff}$ acts on $
\mathfrak{g}$ via
$$[Y, \delta_n] = n \delta_n,  \quad  [X, \delta_n] = \delta_{n+1},$$
for all $n$. Let $\mathfrak{g}_{CM}:=\mathfrak{g}_{aff}\rtimes \mathfrak{g}$ be the
corresponding semidirect product Lie algebra.
As an algebra, $\mathcal{H}_1$ coincides with
the universal enveloping algebra of the Lie algebra $\mathfrak{g}_{CM}$. Thus $\mathcal{H}_1$
is the universal algebra generated by  $\{X, Y, \delta_n ; n =
1, 2, \dots \}$ subject to the relations
$$[Y,X] = X, \; \;[Y, \delta_n] = n \delta_n, \; \;[X, \delta_n] =
\delta_{n+1},\; \;
[\delta_k, \delta_l] = 0,$$
 for $ n, k, l = 1 ,
2, \dots.$
We let the counit of $\mathcal{H}_1$ coincide with the counit of $U(\mathfrak{g}_{CM})$. Its
coproduct
and antipode, however, will be certain deformations of the coproduct and antipode of
$U(\mathfrak{g}_{CM})$ as follows.  Using the universal property of  $U(\mathfrak{g}_{CM})$,
one checks that
the relations
  $$\Delta Y = Y \otimes 1 + 1 \otimes  Y, \quad  \Delta \delta_1 =
  \delta_1 \otimes  1 + 1 \otimes  \delta_1, $$
$$\Delta X = X \otimes  1 + 1 \otimes  X + \delta_1 \otimes  Y,$$
determine a unique  algebra map $ \Delta : \mathcal{H}_1 \to
\mathcal{H}_1 \otimes  \mathcal{H}_1$. Note that $\Delta$ is not
cocommutative and it differs from the coproduct of the enveloping
algebra $U(\mathfrak{g}_{CM})$. Similarly, one checks that there
is a unique antialgebra map, the antipode,  $S: \mathcal{H}_1 \to \mathcal{H}_1$
determined by the relations
$$ S(Y ) = -Y, \; \;  S(X) = -X + \delta_1Y, \; \; S(\delta_1) = -\delta_1.$$

The first realization of $\mathcal{H}_1$ was through its action as
quantum symmetries of the `frame bundle' of the noncommutative space of leaves of codimension one foliations.
More precisely, given a codimension one foliation $(V, \mathcal{F})$, let $M$ be a smooth transversal for
$(V, \mathcal{F})$. Let  $A=C^{\infty}_0(F^+M)$ denote the algebra of smooth
 functions with compact support on the bundle of positively oriented
 frames on $M$ and let   $\Gamma \subset Diff^+(M)$ denote the holonomy group of
  $(V, \mathcal{F})$. One has a natural
 prolongation of the action of  $\Gamma$ to
  $F^+(M)$ by
  $$\varphi (y, y_1)=(\varphi (y), \varphi' (y)(y_1)).$$
 Let $A_{\Gamma}= C^{\infty}_0(F^+M)\rtimes \Gamma$ denote the corresponding crossed product
algebra. Thus the elements of $A_{\Gamma}$ consist of finite linear combinations
(over $\mathbb{C}$)
of terms $fU_{\varphi}^{\ast}$ with $f\in C^{\infty}_0(F^+M)$ and  $\varphi \in \Gamma$. Its product
is defined by
$$fU_{\varphi}^{\ast} \cdot gU_{\psi}^{\ast}=(f \cdot \varphi (g))U_{\psi \varphi}^{\ast}.$$
 There is an
action of  $\mathcal{H}_1$ on $A_{\Gamma}$,   given by \cite{como2, como5}:
$$Y(fU_{\varphi}^{\ast})= y_1 \frac{\partial f}{\partial y_1}
U_{\varphi}^{\ast}, \quad X(fU_{\varphi}^{\ast})= y_1 \frac{\partial f}{\partial y}
U_{\varphi}^{\ast},$$
$$\delta_n(fU_{\varphi}^{\ast})=y_1^n
\frac{ d^n} {dy^n} (log \frac{d\varphi}{dy})fU_{\varphi}^{\ast}.$$
Once these formulas are given, it can be checked, by a long
computation, that $A_{\Gamma}$ is indeed an $\mathcal{H}_1$-module
algebra. To define the corresponding characteristic  map, Connes and Moscovici defined   a modular pair in involution
$(\delta, 1)$ on $\mathcal{H}_1$
and a $\delta$-invariant trace on  $A_{\Gamma}$ as we shall describe next.

 Let $\delta $
denote the unique extension of the modular character
$$\delta : \mathfrak{g}_{aff} \to \mathbb{R}, \quad \delta(X)=1, \; \delta(Y)=0,$$
to a character $\delta : U (\mathfrak{g}_{aff})\to \mathbb{C}.$
There is a unique extension of $\delta $ to a character, denoted
by the same symbol  $\delta : \mathcal{H}_1 \to \mathbb{C}.$
Indeed,  the relations $[Y, \delta_n]=n\delta_n$ show that we must
have $\delta (\delta_n)=0$, for $n=1, 2, \dots.$ One can then
check that these relations are compatible with the algebra
structure of $\mathcal{H}_1.$ The algebra
$A_{\Gamma}=C^{\infty}_0(F^+(M)\rtimes \Gamma$
 admits a $\delta$-invariant trace $\tau : A_{\Gamma} \to
\mathbb{C}$  given by
\cite{como2}:
$$\tau (fU^{\ast}_{\varphi})=\int_{F^+(M)}f(y, y_1)\frac{dy dy_1}{y_1^2},
\quad  \text{if} \; \varphi =1,$$
and $\tau (fU^{\ast}_{\varphi})=0$, otherwise. Now,  using the
$\delta$-invariant trace $\tau$ and the above defined action
$\mathcal{H}_1 \otimes A_{\Gamma}\to A_{\Gamma}$, the characteristic map
\eqref{cmcmap} takes the form
$$ \chi_{\tau}:  HC^*_{(\delta, 1)} (\mathcal{H}_1) \longrightarrow  HC^* (A_{\Gamma}).$$
This map plays a fundamental role in transverse index theory in
\cite{como2}.

The Hopf algebra $\mathcal{H}_1$ shows its beautiful head in number
theory as well. To give an indication of this,
I shall briefly discuss the  {\it modular Hecke algebras} and actions of
$\mathcal{H}_1$ on them
as they were introduced by Connes and Moscovici in \cite{como5, como6}.
 For each
$N\geq 1$, let
$$\Gamma (N)= \left\{
\left( \begin{matrix} a & b\\c & d
\end{matrix} \right) \in SL(2, \mathbb{Z}); \; \; \left( \begin{matrix} a & b\\c & d
\end{matrix}\right)=\left( \begin{matrix} 1 & 0\\0 & 1
\end{matrix}\right) \;  \text{mod}\, N \right\}$$
denote the level $N$ congruence subgroup of $\Gamma (1)=SL(2,
\mathbb{Z})$. Let $\mathcal{M}_k (\Gamma (N))$ denote the space of modular
forms of level $N$ and weight $k$ and
$$\mathcal{M} (\Gamma (N)):=\bigoplus_k \mathcal{M}_k (\Gamma (N))$$
be the graded algebra of modular forms of level $N$. Finally,  let
$$\mathcal{M}:= \underset{\underset{N}{\to}}{\text{lim}} \; \mathcal{M} (\Gamma (N))$$
denote the algebra of modular forms of all levels, where the inductive system is defined by divisibility.
The group
$$G^{+}(\mathbb{Q}):= GL^+ (2, \mathbb{Q}),$$
acts  on  $\mathcal{M}$ through  its usual action on functions on the upper half-plane
(with corresponding weight):
$$(f, \alpha) \mapsto f|_k \alpha (z)= \text{det}(\alpha)^{k/2}(cz+d)^{-k}f(\alpha \cdot
z),$$
$$\alpha = \left( \begin{matrix} a & b\\c & d
\end{matrix}\right), \quad \alpha \cdot z=\frac{az+b}{cz+d}.$$
The simplest modular Hecke algebra is the crossed-product algebra
$$\mathcal{A}=\mathcal{A}_{G^+(\mathbb{Q})}:= \mathcal{M}\rtimes G^+(\mathbb{Q}).$$
Elements of this (noncommutative) algebra will be denoted by  finite sums
$\sum fU^*_{\gamma},$\\ $f\in \mathcal{M}, \, \gamma \in
G^+(\mathbb{Q}).$
$\mathcal{A}$ can be thought of as the algebra of noncommutative
coordinates  on the noncommutative quotient space of modular
forms modulo Hecke correspondences.

Now consider the operator $X$ of degree two on the space of modular forms
defined by
$$X:= \frac{1}{2\pi i}\frac{d}{dz}-\frac{1}{12\pi i}\frac{d}{dz}
(\log
\Delta)\cdot Y,$$
where
$$\Delta (z)=(2\pi)^{12}\eta^{24}(z)=(2\pi)^{12}q\prod_{n=1}^{\infty} (1-q^n)^{24},
\quad q=e^{2\pi iz},$$
$\eta$ is the Dedekind eta function, and $Y$ is  the grading operator
$$Y(f)= \frac{k}{2}\cdot f, \quad \text{for all}\; f\in \mathcal{M}_k.$$
It is shown in \cite{como5}
that there is a unique action of $\mathcal{H}_1$ on
$\mathcal{A}_{G^+(\mathbb{Q})}$ determined by
$$ X(fU^*_{\gamma})=X(f)U^*_{\gamma}, \quad
Y(fU^*_{\gamma})=Y(f)U^*_{\gamma},$$
$$\delta_1( fU^*_{\gamma})=\mu_{\gamma}\cdot f(U^*_{\gamma}),$$
where
$$\mu_{\gamma} (z)=\frac{1}{2\pi
i}\frac{d}{dz} \text{log} \frac{\Delta|\gamma}{\Delta}.$$
This action is compatible with the algebra structure, i.e.,
 $\mathcal{A}_{G^+(\mathbb{Q})}$
is an $\mathcal{H}_1$-module algebra. Thus one can think of
$\mathcal{H}_1$ as quantum symmetries of the noncommutative space
represented by $\mathcal{A}_{G^+(\mathbb{Q})}$.

More generally, for any congruence subgroup $\Gamma$,   an algebra
$A(\Gamma)$ is constructed in \cite{como5} that contains  as
subalgebras both the algebra of $\Gamma$-modular forms and the
Hecke ring at level $\Gamma$. There is also a corresponding action
of $\mathcal{H}_1$ on $A(\Gamma)$.

The Hopf cyclic cohomology groups $HC^n_{(\delta, \sigma)}(H)$ are
computed in several cases in \cite{como2}. Of particular interest
for applications to transverse index theory and number theory is
the (periodic)  cyclic cohomology of $\mathcal{H}_1$. It is shown in \cite{como2} that the
periodic groups $HP^n_{(\delta, 1)}(\mathcal{H}_1)$ are canonically
isomorphic to the Gelfand-Fuchs cohomology, with trivial coefficients,  of the Lie algebra
$\mathfrak{a}_1$ of
formal vector fields on the line:
$$HP^*_{(\delta,
1)}(\mathcal{H}_1)= \bigoplus_{i\geq 0} H^{* +2i}_{GF}(\mathfrak{a}_1, \mathbb{C}). $$
 This result is very significant in that it relates the Gelfand-Fuchs construction of
 characteristic classes of smooth manifolds to a  noncommutative geometric construction of the same
  via
 $\mathcal{H}_1$. Connes and Moscovici also identified certain interesting elements in the Hopf cyclic cohomology of $\mathcal{H}_1$.
   For example,
it can be directly checked that the elements $\delta_2' :
=\delta_2-\frac{1}{2}\delta_1^2$ and $\delta_1$ are Hopf cyclic 1-cocycles
for  $\mathcal{H}_1$, and
$$F:=X\otimes Y-Y\otimes X -\delta_1Y\otimes Y$$
is a  Hopf cyclic 2-cocycle. Under the characteristic map \eqref{cmcmap} and for $A= A_{\Gamma}$ these Hopf
cyclic cocycles are mapped to  the Schwarzian derivative, the Godbillon-Vey cocycle, and
 the transverse fundamental class of Connes \cite{ac83c}, respectively. See
 \cite{como6} for detailed calculations  as well as  relations with
 modular forms and modular Hecke algebras.
Very recently the  unstable  cyclic cohomology groups of $\mathcal{H}_1$, and a series of
other Hopf algebras
attached to pseudogroups of geometric structures,   were fully
computed in \cite{moscrang1, moscrang2}. In particular it is shown that
the groups $HC^n_{(\delta, \sigma)}(\mathcal{H}_1)$ are finite
dimensional for all $n$.

The notion of  modular pair in involution  $(\delta, \sigma)$ for a Hopf algebra might seem rather ad hoc at a first glance.
This is in fact not the case and the concept is very natural and fundamental.
For example,  it is shown in \cite{como3} that {\it coribbon Hopf algebras} and compact
quantum groups are endowed with canonical modular pairs of the form $(\delta, 1)$
and,
dually,  {\it ribbon Hopf algebras}  have canonical modular pairs of the type
$(1, \sigma)$. The fundamental importance of modular pairs in involution was further elucidated when
{\it Hopf cyclic cohomology  with coefficients} was   introduced in \cite{hkrs1, hkrs2}. It turns out that
some  very stringent  conditions
have to be imposed on an $H$-module $M$ in order for $M$ to serve as a
coefficient (local system)  for  Hopf cyclic cohomology theory. Such
modules are called {\it stable anti-Yetter-Drinfeld modules}. More precisely,
a (left-left)  anti-Yetter-Drinfeld  $H$-module is a left $H$-module $M$
which is simultaneously a
left $H$-comodule such that
$$\rho(hm)=h^{(1)}m^{(-1)}S(h^{(3)})\ot h^{(2)}m^{(0)},$$
for all $h\in H$ and  $m\in M.$ Here $\rho: M \to H \otimes M$, $\rho(m)= m^{(-1)}\otimes m^{(0)}$ is  the
comodule structure map  of $M$.
 $M$ is called stable if in addition we have
$$m^{(-1)}m^{(0)}=m,  $$
for all $m\in M$.
Given a stable  anti-Yetter-Drinfeld (SAYD) module  $M$ over $H$, one  can then define the  Hopf cyclic cohomology of
$H$ with coefficients in M. 
One-dimensional
SAYD modules correspond to Connes-Moscovici's modular pairs in
involution. More precisely,
there is a one-to-one correspondence between modular pairs in
involution $(\delta, \sigma)$ on $H$ and SAYD module structures on $M=\mathbb{C}$, the ground field, defined by
$$h. r=\delta (h) r, \quad r\mapsto \sigma \otimes r,$$
for all $h\in H$ and $r\in \mathbb{C}$. Thus a modular pair in
involution can be regarded as an `equivariant line bundle' over the
noncommutative space represented by the Hopf algebra $H$.

{}


\begin{thebibliography}{999}			




\bibitem{atising}  M. F. Atiyah and I. M.  Singer,   The index of elliptic operators. IV.
Ann. of Math. (2),  93  1971,  119-138.












\bibitem{bur}
 D. Burghelea, The cyclic homology of the group rings.
 Comment. Math. Helv.  {\bf 60}  (1985),  no. 3, 354--365.











\bibitem{ac78}
 A. Connes, A survey of foliations and operator algebras. Operator algebras and applications,
 Part I (Kingston, Ont., 1980), pp. 521--628, Proc. Sympos. Pure Math., 38, Amer. Math. Soc.,
 Providence, R.I., 1982.



\bibitem{ac80}
A. Connes, { $C\sp{*} $ alg\'ebres et g\'eom\'etrie diff\'erentielle}. (French) C. R. Acad.
 Sci. Paris S\'er. A-B  {\bf 290} (1980), no. 13, A599--A604.

\bibitem{ac81}
A. Connes, {Spectral sequence and homology of currents for operator
algebras}. Oberwolfach Tagungsber., 41/81, Funktionalanalysis und
$C^*$-Algebren, 27-9/3-10, 1981.


\bibitem{ac83a}
 A. Connes, {Noncommutative differential geometry}. Chapter I: The Chern character in K homology,
 Preprint IHES octobre 1982 ; Chapter II: de Rham homology and
 noncommutative algebra, Preprint IHES f\'{e}vrier 1983.  

\bibitem{ac83b}
A. Connes, {Cohomologie cyclique et foncteurs  ${\rm Ext}\sp n$ }. (French) (Cyclic cohomology
and functors ${\rm Ext}\sp n$) C. R.
 Acad. Sci. Paris Sér. I Math. {\bf 296} (1983), no. 23, 953--958.

\bibitem{ac83c}
 A. Connes,  Cyclic cohomology and the transverse fundamental class of a foliation.
 Geometric methods in operator algebras (Kyoto, 1983),  52--144, Pitman Res. Notes Math. Ser., 123,
 Longman Sci. Tech., Harlow, 1986.


\bibitem{ac85}
 A. Connes, {Noncommutative differential geometry}. Chapter I: The Chern character in K homology;
  Chapter II: de Rham homology and
 noncommutative algebra,  Publ. Math. IHES no.
 62 (1985), 41-144.



\bibitem{ac88}A. Connes, Entire cyclic cohomology of Banach algebras and characters
of $\theta$-summable Fredholm modules.  $K$-Theory  1  (1988),  no. 6, 519--548.


\bibitem{acaction}  A. Connes, The action functional in noncommutative geometry.
Comm. Math. Phys.  117  (1988),  no. 4, 673--683.

\bibitem{acmk}
A. Connes and M.  Karoubi,  Caract\`{e}re multiplicatif d'un module de Fredholm. (French)
[Multiplicative character of a Fredholm module]  $K$-Theory  2  (1988),  no. 3, 431--463.



\bibitem{acbook}
 A. Connes, {Noncommutative geometry}. Academic Press, Inc., San Diego, CA, 1994.

\bibitem{ac99}
A. Connes,  Trace formula in noncommutative geometry and the zeros of the Riemann zeta function.
Selecta Math. (N.S.)  5  (1999),  no. 1, 29--106.


\bibitem{ac00}
A. Connes,  Noncommutative geometry year 2000.
Highlights of mathematical physics (London, 2000),  49--110, Amer. Math.
Soc., Providence, RI, 2002.

\bibitem{cocoma}
A. Connes, C. Consani  and M. Marcolli, Noncommutative geometry and motives:
the thermodynamics of endomotives, Advances in Math. 214 (2) (2007).











 \bibitem{coma3} A. Connes  and M. Marcolli,  Noncommutative geometry,
 quantum fields and motives. American Mathematical Society Colloquium Publications, 55. American Mathematical
 Society, Providence, RI; Hindustan Book Agency, New Delhi, 2008.
 xxii+785 pp.



\bibitem{como0}
A. Connes and H. Moscovici, Cyclic cohomology, the Novikov conjecture and hyperbolic groups.
Topology  29  (1990),  no. 3, 345--388.


\bibitem{como1}
 A. Connes   and H. Moscovici,  The local index formula in noncommutative geometry.  Geom. Funct.
 Anal.  {\bf 5}  (1995),  no. 2, 174--243.

\bibitem{como2} A. Connes   and H. Moscovici, {Hopf algebras, cyclic cohomology and the transverse
index theorem}, Comm. Math. Phys. {\bf 198} (1998), no. 1, 199--246.

\bibitem{como3}
A. Connes   and   H. Moscovici, {Cyclic cohomology and Hopf algebra symmetry.}
Conference Mosh\'{e} Flato 1999 (Dijon).
 Lett.
Math. Phys. {\bf 52} (2000), no. 1, 1--28.

\bibitem{como4}
A. Connes   and   H. Moscovici,  Rankin-Cohen brackets and
the Hopf algebra of transverse geometry. Mosc. Math. J. 4 (2004), no. 1,
111--130.

\bibitem{como5}
A. Connes   and   H. Moscovici,   Modular Hecke
algebras and their Hopf symmetry. Mosc. Math. J. 4 (2004), no. 1,
67--109, 310.

\bibitem{como6} A. Connes   and H. Moscovici, Background independent
geometry and Hopf cyclic cohomology, math.QA/0505475.

\bibitem{como7} A. Connes   and H. Moscovici, Type III and spectral triples.
 Traces in number theory, geometry and quantum fields,  57--71, Aspects Math., E38, Friedr. Vieweg, Wiesbaden, 2008.







\bibitem{contre}
A. Connes and P. Tretkoff, The Gauss-Bonnet theorem for the
noncommutative two torus, arXiv:0910.0188v1.

\bibitem{fakh1} F. Fathi Zadeh and M. Khalkhali, The
algebra of formal twisted pseudodifferential symbols and a noncommutative residue, arXiv:0810.0484. 

\bibitem{fakh2} F. Fathi Zadeh and M. Khalkhali, The
Gauss-Bonnet theorem for noncommutative two tori with a general
conformal structure,  arXiv:1005.4947. To appear in Journal of
Noncommutative Geometry. 



\bibitem{freed} M. Freedman, Towards quantum characteristic classes, Lecture at MSRI, August 11, 2008,
available at
http://www.msri.org/communications/vmath/VMathVideos/VideoInfo/3895/-\\showvideo.




\bibitem{hkrs1}
P. Hajac, M. Khalkhali, B. Rangipour, and Y. Sommerh\"{a}user,
Stable anti-Yetter-Drinfeld modules.  C. R. Math. Acad. Sci. Paris  338  (2004),  no. 8, 587--590.

\bibitem{hkrs2}
P. Hajac, M. Khalkhali, B. Rangipour, and Y. Sommerh\"{a}user,
Hopf-cyclic homology and cohomology with coefficients.  C. R. Math. Acad. Sci. Paris  338
(2004),  no. 9, 667--672.



\bibitem{kh}
M. Khalkhali, Basic Noncommutative Geometry, European Mathematical
Society Publishing House, 2009. 






\bibitem{lod}
     J. L. Loday, {Cyclic Homology.} Springer-Verlag, second edition (1998).

\bibitem{mosc1} H. Moscovici, Local index formula and twisted spectral triples, arXiv:0902.0835. This volume. 

\bibitem{moscrang1} H. Moscovici,   and B. Rangipour, Cyclic cohomology of Hopf algebras of transverse
symmetries in codimension 1.  Adv. Math.  210  (2007),  no. 1, 323--374. 

\bibitem{moscrang2} H. Moscovici,  and B. Rangipour Hopf algebras of primitive
Lie pseudogroups and Hopf cyclic cohomology.  Adv. Math.  220  (2009),
no. 3, 706--790.


\end{thebibliography}
\end{document}